\documentclass{article}

     \PassOptionsToPackage{numbers, compress}{natbib}
 \usepackage[final]{neurips_2019}

\usepackage{algorithmicx,algpseudocode}
\usepackage[utf8]{inputenc} %
\usepackage[T1]{fontenc}    %
\usepackage[colorlinks,citecolor=blue,urlcolor=black,linkcolor=blue,pdfborder={0 0 0}]{hyperref}
\usepackage{url}            %
\usepackage{booktabs}       %
\usepackage{amsfonts}       %
\usepackage{nicefrac}       %
\usepackage{microtype}      %
\usepackage{fullpage, url,amsmath,amsfonts,amssymb,mathtools,mathrsfs,graphicx, algorithm, float, sansmath,epstopdf,xcolor,caption,enumitem,tabularx,subcaption}
\usepackage[capitalize]{cleveref}  
\crefname{nlem}{Lemma}{Lemmas}
\crefname{nprop}{Proposition}{Propositions}
\crefname{ncor}{Corollary}{Corollaries}
\crefname{nthm}{Theorem}{Theorems}
\crefname{exa}{Example}{Examples}
\crefname{assumption}{Assumption}{Assumptions}
\crefname{equation}{}{}
\crefname{enumi}{}{}

\newcommand{\R}{{\mathbb R}}
\newcommand{\E}{{\mathcal E}}

\usepackage{bbm}

 \newcommand{\qtext}[1]{\quad\text{#1}\quad} 
 \newenvironment{talign*}
 {\csname align*\endcsname}
 {\endalign}
\newenvironment{talign}
 {\csname align\endcsname}
 {\endalign}

\newtheorem{definition}{Definition}
\newtheorem{remark}{Remark}

\newtheorem{example}{Example}
\newtheorem{theorem}{Theorem}
\newtheorem{lemma}[theorem]{Lemma}

\newtheorem{proposition}[theorem]{Proposition}

\newcommand{\X}{{\mathcal X}}
\title{Accelerating Rescaled Gradient Descent: \\Fast Optimization of Smooth Functions}

\author{%
  Ashia C.~Wilson\\%\thanks{Use footnote for providing further information
  Microsoft Research \\
  \texttt{ashia.wilson@microsoft.com} \\
   \And
   Lester Mackey \\
    Microsoft Research \\
   \texttt{lmackey@microsoft.com} \\
   \And
   Andre Wibisono \\
   Georgia Tech \\
   \texttt{wibisono@gatech.edu}
}

\begin{document}

\maketitle

\begin{abstract}
We present a family of algorithms, called {\em descent algorithms}, for optimizing convex and non-convex functions. We also introduce a new first-order algorithm, called {\em rescaled gradient descent }(RGD), and show that RGD achieves a faster convergence rate than gradient descent over the class of {\em strongly smooth} functions -- a natural generalization of the standard smoothness assumption on the objective function. When the objective function is convex, we present two frameworks for accelerating descent algorithms, one in the style of Nesterov and the other in the style of Monteiro and Svaiter, using a single Lyapunov function. Rescaled gradient descent can be accelerated under the same strong smoothness assumption using both frameworks. We provide several examples of strongly smooth loss functions in machine learning and numerical experiments that verify our theoretical findings. We also present several extensions of our novel Lyapunov framework including deriving optimal universal higher-order tensor methods and extending our framework to the coordinate descent setting. 

\end{abstract}

\section{Introduction}
We consider the optimization problem 
\begin{align}\label{eq:obj}
\min_{x \in \mathcal{X}} f(x)
\end{align}
where $f: \mathcal{X} \rightarrow \mathbb{R}$ is a continuously differentiable function, %
on  a finite-dimensional real vector space $\mathcal{X}$ with inner product norm $\| v \| :=\sqrt{\langle v,Bv\rangle}$ and a dual norm $\|s \|_{\ast} := \sqrt{\langle s, B^{-1}s\rangle}$ for $s$ in the dual space $\X^\ast$.
Here, $B: \mathcal{X} \rightarrow \mathcal{X}^\ast$ is a positive definite self-adjoint operator. We assume the minimum of $f$ is attainable and let $x^\ast$ represent a point in $\arg \min_{x \in \X} f(x)$. 

We study the performance 
of a family of discrete-time algorithms parameterized by  $\delta >0$ and an integer scalar $ 1<p \leq \infty$, called {\em $\delta$-descent algorithms of order $p$}. These algorithms meet a progress condition that allows us to derive fast non-asymptotic convergence rate upper bounds, parameterized by $p$, for both nonconvex and convex instances of~\cref{eq:obj}. For example, descent algorithms of order $1<p<\infty$ satisfy the upper bound $f(x_k) - f(x^\ast)=O( 1/(\delta k)^{p-1})$ for convex functions.

Using this framework we introduce a new method for smooth optimization called {\em rescaled gradient descent} (RGD),
\begin{align*} 
x_{k+1} = x_k - \eta^{\frac{1}{p-1}} \frac{B^{-1}\nabla f(x_k)}{\|\nabla f(x_k)\|_{\ast}^{\frac{p-2}{p-1}}},\qquad \eta >0, p>1.
\end{align*} 
 We show that if~\cref{eq:obj} is sufficiently smooth, rescaled gradient descent is a $\delta$-descent algorithm of order $p$, and subsequently converges quickly to solutions of~\cref{eq:obj}. RGD can be viewed as a natural generalization of gradient descent ($p=2$) and normalized gradient descent ($p=\infty$), whose non-asymptotic behavior for quasi-convex functions has been well-studied~(\cite{Hazan}). 

When $f$ is convex, we present two frameworks for obtaining algorithms with faster convergence rate upper bounds. The first, pioneered in~\citet{Nesterov83,Nesterov04,Nesterov05,Nesterov08}, shows how to wrap a $\delta$-descent method of order $1<p<\infty$ in two sequences to obtain a method that satisfies $f(x_k) - f(x^\ast) =O( 1/(\delta k)^p)$. The second, introduced by~\cite{Svaiter}, shows how to wrap a $\delta$-descent method of order $1<p<\infty$ in the same set of sequences and add a line search step to obtain a method that satisfies  $f(x_k) - f(x^\ast) =O( 1/(\delta k)^{\frac{3p-2}{2}})$. We provide a general description of both frameworks and show how they can be applied to RGD and other descent methods of order $p$.

Our motivation also comes from a burgeoning literature (e.g.,~\cite{Polyak1964,Schropp,SuBoydCandes14,Accel2,KricheneStoch,Wilson,BetancourtJordanWilson,jelena,Andre,SuJordan,PDE,Hamiltonian}) that harnesses the connection between dynamical systems and optimization algorithms to develop new analyses and optimization methods. 
 Rescaled gradient descent is obtained by discretizing an ODE called {\em  rescaled gradient flow} introduced by~\cite{Acceleration}. 
We compare RGD and accelerated RGD to the work of \citet{Runge},  who introduce accelerated dynamics and apply Runge-Kutta integrators to discretize them. They show that Runge-Kutta integrators converge quickly when the function is sufficiently smooth and when the order of the integrator is sufficiently large. We provide a better convergence rate upper bound for accelerated RGD under a very similar smoothness assumption. We also compare our work to~\citet{Hamiltonian}, who introduces conformal Hamiltonian dynamics and show that if the objective function is sufficiently smooth, algorithms obtained by discretizing these dynamics converge at a linear rate. We show (accelerated) RGD also achieves a fast linear rate under similar smoothness conditions.

The remainder of this paper is organized as follows. %
\cref{sec:2} introduces $\delta$-descent algorithms and  \cref{sec:21} describes several examples of descent algorithms that are popular in optimization. \cref{sec:22} introduces RGD and \cref{sec:3} presents two frameworks for accelerating $\delta$-descent methods and applies both to RGD. \cref{sec:exam} describes several examples of strongly smooth objective functions as well as experiments to verify our findings.  Finally, \cref{sec:4} discusses simple extensions of our framework, including deriving and analyzing {\em optimal universal tensor methods} for objective functions that have H\"older-continuous higher-order gradients and extending our entire framework and results to the coordinate setting.

\section{Descent Algorithms}
\label{sec:2}
The focus of this section is a family of algorithms called $\delta$-{\em descent algorithms of order p}.
\begin{definition}\label{def:descent-alg} An algorithm $x_{k+1} = \mathcal{A}(x_k)$ is a $\delta$-{\bf descent algorithm of order $ p$} for $1<p\leq\infty$ if for some constant $0<\delta<\infty$ it satisfies
\begin{subequations}\label{eq:prog1}
\begin{talign}\label{eq:p1}
\frac{f(x_{k+1}) - f(x_k)}{\delta}&\leq -\|\nabla f(x_{k})\|_{\ast}^{\frac{p}{p-1}} \qtext{for all $k\geq 0$,  or } \\
\label{eq:p2}
\frac{f(x_{k+1}) - f(x_k)}{\delta}&\leq - \|\nabla f(x_{k+1})\|_{\ast}^{\frac{p}{p-1}} \qtext{for all $k \geq 0$. }
\end{talign}
\end{subequations}
\end{definition}
For $\delta$-descent algorithms of order $p$, it is possible to obtain non-asymptotic convergence guarantees for non-convex, convex and gradient dominated functions.
 Recall, a function is \emph{$\mu$-gradient dominated of order $p \in (1,\infty]$} if
\begin{talign}\label{eq:gdom}
\frac{p-1}{p} \|\nabla f(x)\|_{\ast}^{\frac{p}{p-1}} \geq  \mu^{\frac{1}{p-1}} (f(x) - f(x^\ast)), \quad \forall x \in \mathcal{X}.
\end{talign}
When $p = 2$,~\eqref{eq:gdom} is the Polyak-\L{}ojasiewicz condition introduced concurrently by~\citet{Polyak1964} and ~\citet{Loj64}. For the following three theorems, we use the shorthand $E_0 := f(x_0) - f(x^\ast)$ and assume $f$ is differentiable. %
\begin{theorem}\label{thm:11}
Any $\delta$-descent algorithm of order $p$ satisfies
\begin{talign}\label{eq:b11}
\min_{0\leq s \leq k} \|\nabla f(x_s)\|_{\ast} \leq (E_0/(\delta k))^{\frac{p-1}{p}}. %
\end{talign}
\end{theorem}
\begin{theorem}\label{thm:22}
If $f$ is convex with $R = \sup_{x:f(x)\leq f(x_0)} \|x - x^\ast \| < \infty$, and $c_p := \frac{(1-1/p)^p}{p-1}$, then any $\delta$-descent algorithm of order $p$ satisfies
\begin{talign}\label{eq:b22}
f(x_k) - f(x^\ast) 
	\leq 
	\begin{cases}
	2\big(\frac{1}{E_0^{1/p}}+\frac{1}{R\gamma c_p^{1/p}p}(\delta k)^{\frac{p-1}{p}}\big)^{-p}
	=\textstyle O(1/(1+\frac{1}{R\gamma p}{(\delta k)^{\frac{p-1}{p}}})^{p}), & p < \infty \\
	2E_0 \exp(-\delta k/(R\gamma)), & p = \infty.
	\end{cases}
\end{talign}
where $\gamma = 1$ when \cref{eq:p1} is satisfied and $\gamma = (1+\frac{1}{Rp}(E_0/c_p)^{
\frac{1}{p}} \delta^\frac{p-1}{p})^{p-1}$ when \cref{eq:p2} is satisfied.

\end{theorem}
\begin{theorem}\label{thm:33}
If $f$ is $\mu$-gradient dominated of order $p$, then
any $\delta$-descent algorithm of order $p$ satisfies
\begin{talign}\label{eq:b33}
f(x_k) - f(x^\ast) \leq E_0\exp\left(-\frac{p}{p-1}\mu^{\frac {1}{p-1}}\delta k\right). %
\end{talign}
\end{theorem}
The proof of \cref{thm:11,thm:22,thm:33} are all based on simple energy arguments and can be found in \cref{App:A}. 
Bounds of the form \cref{eq:b11} are common in the non-convex optimization literature and have previously been established for gradient descent ($p=2$ see e.g.~\cite[Thm1]{Nesterov06}) and higher-order tensor methods (see e.g.\cite{Yair}). \cref{thm:11} provides a more general description of algorithms that satisfy this kind of bound.  

Typically, algorithms satisfy the progress condition~\cref{eq:prog1} for specific smoothness classes of functions. 
For example, gradient descent with step-size $0<\eta \leq 1/L$ {\em is a $\delta$-descent method of order $p=2$} with $\delta= \eta/2$ when $\|\nabla^2 f\|\leq L$.  Throughout, we denote $\|B\| = \max_{\|h\|\leq 1} \|B h\|_\ast$, for any $B: \mathcal{X} \rightarrow \mathcal{X}^\ast$.
We list several other examples. 

\subsection{Examples of descent algorithms}
\label{sec:21}
\cref{thm:11,,thm:22,,thm:33} provide a seamless way to derive standard upper bounds for many algorithms in optimization.
\begin{example} The {\em universal higher-order tensor method},
\begin{align}\label{eq:hgd}
x_{k+1} = \arg \min_{x\in \mathcal{X}} \left\{ f_{p-1}(x; x_k) + \frac{1}{\tilde{p}\eta }\|x - x_k\|^{\tilde p}\right\},
\end{align}
where $f_{p-1}(y; x) = \sum_{i=0}^{p-1} \frac{1}{i!}\nabla^i f(x)(y-x)^i$
is the $(p-1)$-st order Taylor approximation of $f$ centered at $x$ and $\tilde p = p -1 + \nu$ for $\nu \in (0,1]$, has been studied by several works %
\cite{Baes09,Acceleration,Nesterov18}.
When $f$ is convex and has H\"older-smooth $(p-1)$-st order gradients, namely $\|\nabla^{p-1} f(x) -\nabla^{p-1} f(y)\|\leq L\|x-y\|^\nu$,~\cref{eq:hgd} with step size $0 < \eta \leq \frac{\sqrt{3}(p-2)!}{2L}$, is a $\delta$-descent algorithm of order $\tilde p$ with $\delta = \eta^{\frac{1}{\tilde p-1}}/2^{\frac{2\tilde p-3}{\tilde p-1}}$.

\end{example}
\begin{example} {\em The natural proximal method},
\begin{align}\label{eq:prox}
x_{k+1} = \arg \min_{x \in \mathcal{X}}\left\{ f(x) + \frac{1}{p\eta }\|x - x_k\|_{x_k}^p\right\},
\end{align}
 where $\| v\|_{x} = \sqrt{\langle v, \nabla^2 h(x)v\rangle}$ was introduced in the setting $h(x) = \frac{1}{2}\|x\|_2^2$ by~\cite{prox}.
For any $\eta, m > 0$ and $mB\preceq \nabla^2 h $, the proximal method is a $\delta$-descent algorithm of order $p$ with $\delta = m^{\frac{p}{p-1}}\eta^{\frac{1}{p-1}}/p$.
\end{example}

\begin{example} {\em Natural gradient descent},
\begin{align}\label{eq:ngd}
x_{k+1} = x_k - \eta  \nabla^2 h(x_k)^{-1}\nabla f(x_k) 
= \arg \min_{x \in \mathcal{X}} \left\{\langle \nabla f(x_k), x\rangle+\frac{1}{2\eta}\|x - x_k\|_{ x_k}^2\right\},
\end{align}
where $\| v\|_{x} = \sqrt{\langle v, \nabla^2 h(x)v\rangle}$ was introduced by~\cite{Amari1988}. %
Suppose $\|\nabla^2 f\| \leq L$ and $m B \preceq \nabla^2 h \preceq MB$ for some $m, L, M > 0$.
Then natural gradient descent with step size $0 < \eta \le \frac{m^2}{ML}$  is a $\delta$-descent algorithm of order $p=2$ with $\delta = \frac{\eta}{2M}$.
\end{example}

\begin{example} {\em Mirror descent},
\begin{align}\label{eq:md}
x_{k+1} =  \arg \min_{x \in \mathcal{X}} \left\{\langle \nabla f(x_k), x\rangle+\frac{1}{\eta}D_h(x ,x_k)\right\},
\end{align}
 where $D_h(x,y) = h(x) - h(y) - \langle \nabla h(y), x-y\rangle$ is the Bregman divergence was introduced by~\cite{NemirovskiiYudin}. 
Suppose $\|\nabla^2 f\| \leq L$ and $m B \preceq \nabla^2 h \preceq MB$ for some $m, L, M > 0$.
Then mirror descent with step size $0 < \eta \le \frac{m^2}{ML}$  is a $\delta$-descent algorithm of order $p=2$ with $\delta = \frac{\eta}{2M}$.
\end{example}

\begin{example} {\em The proximal Bregman method},
\begin{align}\label{eq:pbm}
x_{k+1} =  \arg \min_{x \in \mathcal{X}} \left\{ f(x)+\frac{1}{\eta}D_h(x ,x_k)\right\},
\end{align}
 was introduced by~\cite{ChenTeboulle93}).
 When $m B \preceq \nabla^2 h \preceq MB$ the proximal Bregman method with step-size $\eta >0$ is a $\delta$-descent algorithm of order $p=2$ with $\delta = \frac{m\eta}{2M^2}$.

\end{example}

Details for these examples are contained in \cref{sec:examples}. %

\subsection{Rescaled gradient descent} 
\label{sec:22}
We end this section by discussing the function class for which 
rescaled gradient descent (RGD), 
\begin{talign}\label{eq:rgd}
x_{k+1}  = x_k - \eta^{\frac{1}{p-1}}\frac{B^{-1}\nabla f(x_k)}{\|\nabla f(x_k)\|_{\ast}^{\frac{p-2}{p-1}}} = \arg \min_{x\in\X}\left\{\langle \nabla f(x_k), x\rangle + \frac{1}{p\eta }\|x - x_k\|^p\right\},
\end{talign}
is a $\delta$-descent method of order $p$.
\begin{definition}\label{ass:ss2}
A function $f$ is {\bf strongly smooth} of order $p$ for some  integer $p > 1$, if there exist constants $0<L_1, \dots, L_p < \infty$ such that for $m = 1,\dots,p-1$ and for all $x \in \mathbb{R}^d$:
\begin{align}\label{eq:ss}
|\nabla^m f(x)(B^{-1}\nabla f(x))^m| \leq L_m \|\nabla f(x)\|_{\ast}^{m + \frac{p-m}{p-1}}
\end{align} 
and moreover for $m=p$, $f$ satisfies the condition
$|\nabla^p f(x)(v)^p| \leq L_p \|v\|^p$, $\forall v \in \mathcal{X}$. %
\end{definition}
Here, $\nabla^m f(x)( h)^m  = \sum_{i_1, \dots, i_m = 1}^d \partial_{x_{i_1} \dots x_{i_m}} f(x) \prod_{j=1}^m h_{i_j}$ where $\partial_{x_i} f$ is the partial derivative of $f$ with respect to $x_i$.
 We can always take $L_1 = 1$. 
When $p=2$,~\cref{eq:ss} is the usual Lipschitz condition on the gradient of $f$, but otherwise~\cref{eq:ss}  is stronger.
In particular, if $f$ is strongly smooth of order $p$, then the minimizer $x^\ast$ has order at least $p-1$, i.e., the higher gradients vanish: $\nabla^m f(x^\ast) = 0$ for $m = 1,\dots,p-1$,
whereas this is not implied under mere smoothness. An example of a strongly smooth function of order $p$ is the $p$-th power of the $\ell_2$-norm 
$f(x) = \|x\|_2^p $ with $B = I$, or the $\ell_p$-norm $f(x) = \|x\|_p^p$.
We discuss other families of strongly smooth functions in \cref{sec:exam}. Finally, it is worth mentioning that for most of our results, the absolute value on the left hand side of~\eqref{eq:ss} is unnecessary. We now present the main result regarding the performance of RGD on functions that satisfy~\cref{eq:ss}:

\begin{theorem}\label{lem:1}
Suppose $f$  is strongly smooth of order $p > 1$ with constants $0<L_1, \dots, L_p < \infty$. Then rescaled gradient descent with step-size
\begin{talign}\label{eq:step}
0 < \eta^{\frac{1}{p-1}} \leq \min\left\{1, \frac{1}{\left(2\sum_{m=2}^p \frac{L_m}{m!}\right)}\right\}
\end{talign}
 satisfies the descent condition~\cref{eq:p1} with $\delta = \eta^{\frac{1}{p-1}}/2$. 
\end{theorem}
The proof of \cref{lem:1} is in \cref{App:RGD}.  A corollary to Theorems~\ref{thm:11}-\ref{lem:1} is the following theorem.
\begin{theorem}
RGD with a step size that satisfies~\cref{eq:step} achieves convergence rate guarantee~\cref{eq:b11} when $f$ is differentiable and strongly smooth of order p,~\cref{eq:b22} when $f$ is convex function and strongly smooth  of order p, and~\cref{eq:b33} when $f$ is $\mu$-uniformly convex and strongly smooth of order p, where $\delta^{p-1} = \eta/2^{p-1}$.
\end{theorem}  

Our results show rescaled gradient descent can minimize the canonical $p$-strongly smooth and uniformly convex function $f(x) = \frac{1}{p} \|x\|^p$ at an exponential rate;
in contrast, gradient descent can only minimize it at a polynomial rate, even in one dimension.
We provide the proof of Proposition~\ref{Prop:GDSlow} in Appendix~\ref{App:GDSlow}.

\begin{proposition}\label{Prop:GDSlow}
Let $f \colon \R \to \R$ be $f(x) = \frac{1}{p} |x|^p$ for $p > 2$, with minimizer $x^\ast = 0$ and $f(x^\ast) = 0$.
For any step size $0 < \eta^{\frac{1}{p-1}} < 1$ and initial position $x_0 \in \R$, rescaled gradient descent of order $p$ minimizes $f$ at an exponential rate: $f(x_k) = (1-\eta^{\frac{1}{p-1}})^{pk} f(x_0)$.
On the other hand, for  any $\eta^{\frac{1}{p-1}} > 0$ and $|x_0| < (2\eta^{\frac{1}{p-1}})^{-\frac{1}{p-2}}$, gradient descent minimizes $f$ at a polynomial rate: $f(x_k) = \Omega((\eta^{\frac{1}{p-1}} k)^{-\frac{p}{p-2}})$.
\end{proposition}

 We now demonstrate how all the aforementioned examples of $\delta$-descent methods can be accelerated.

 \section{Accelerating Descent Algorithms}
 \label{sec:3}
We present two frameworks for accelerating descent algorithms based on the dynamical systems perspective introduced by~\citet{Acceleration} and~\citet{Wilson} and apply them to RGD. The backbone of both frameworks is the Lyapunov function 
\begin{talign*}
 E_k = A_k(f(x_k) - f(x^\ast)) + D_h(x^\ast, z_k),
 \end{talign*} and two sequences~\eqref{eq:nest1} and~\eqref{eq:nest2}. The connection between continuous time dynamical systems and these two sequences and Lyapunov function is described in~\cite{Wilson}. We present a high-level description of both techniques in the main text and leave details %
 of our analysis to~\cref{app:accel}. %
\subsection{Nesterov acceleration of descent algorithms}

In the context of convex optimization, the technique of ``acceleration''  has its origins in~\citet{Nesterov83} and refined in~\citet{Nesterov04}.
In these works, Nesterov showed how to combine gradient descent with two sequences 
 to obtain an algorithm with an optimal convergence rate. There have been many works since (as well as some frameworks, including~\cite{Zaid,coupling,Lessard14,Wilson}) describing how to accelerate various other algorithms to obtain methods with superior convergence rates. 

 \citet{Wilson}, for example, show the following two discretizing schemes,
 \begin{subequations}\label{eq:nest1}
\begin{talign}
x_{k} &= \delta\tau_k z_k + (1 - \delta\tau_k) y_k\label{eq:up1}\\
z_{k+1} &= \arg \min_{z \in \mathcal{X}}\left\{ \alpha_k\langle \nabla f(x_k), z\rangle +  \frac{1}{\delta} D_h(z, z_k)\right\}\label{eq:up2}
\end{talign}
\end{subequations} 
where $y_{k+1}$  satisfies the $\delta^{\frac{p}{p-1}}$-descent condition  $f(y_{k+1}) - f(x_{k}) \leq -  \delta^{\frac{p}{p-1}}\|\nabla f(x_{k})\|_{\ast}^{\frac{p}{p-1}};$ and 
\begin{subequations}\label{eq:nest2}
\begin{talign}
x_{k} &= \delta\tau_k z_k + (1 - \delta\tau_k) y_k\label{eq:up11}\\
z_{k+1} &= \arg \min_{z}\left\{ \alpha_k \langle \nabla f(y_{k+1}), z\rangle + \frac{1}{\delta} D_h(z, z_k)\right\}\label{eq:up22},
\end{talign}
\end{subequations}
where the update for $y_{k+1}$ satisfies the ($\delta^{\frac{p}{p-1}}$-descent) condition %
$f(y_{k+1}) - f(x_{k})\leq  \langle \nabla f(y_{k+1}), y_{k+1} - x_k\rangle \leq - \delta^{\frac{p}{p-1}} \|\nabla f(y_{k+1})\|_{\ast}^{\frac{p}{p-1}},$ constitute an ``accelerated method''.
Their results can be summarized in the following theorem.
 \begin{theorem}\label{prop:8} Assume for all $x,y\in \mathcal{X}$, the function $h$ satisfies the local uniform convexity condition $D_h(x,y) \geq \frac{1}{p}\|x - y\|^p$. Then sequences~\cref{eq:nest1} and~\cref{eq:nest2} with parameter choices $\alpha_k = (\delta/p)^{p-1} k^{(p-1)}$ (where $ k^{(p)} := k(k+1)\cdots(k+p-1)$ is the rising factorial) and $\tau_k = \frac{p}{\delta(p+k)} = \Theta(\frac{p}{\delta k})$ satisfy, 
 \begin{talign}\label{eq:rate1}
 f(y_{k}) - f(x^\ast) \leq \frac{p^pD_h(x^\ast, \,z_0)}{(\delta k)^p} = O\left(1/(\delta k)^p\right).
 \end{talign}
 \end{theorem}
Proof details are contained in \cref{app:acc}. \citet{Wilson} call these new methods {\em accelerated descent methods} due to the fact that \cref{thm:22} guarantees implementing just the $y_{k+1}$ sequence (where we set $x_{k} = y_{k}$) satisfies $f(y_k) - f(x^\ast) \leq O(1/(\tilde\delta k)^{p-1})$,
 where $\tilde \delta^{p-1} = \delta^{p}$. The computational cost of adding sequences~\cref{eq:up1} and~\cref{eq:up2} (or~\cref{eq:up11} and~\cref{eq:up22}) to the descent method is at most an additional gradient evaluation. 
 \begin{remark}[Restarting for accelerated linear convergence]If, in addition, $f$ is $\mu$-gradient dominated of order $p$, then algorithms~\cref{eq:nest1} and~\cref{eq:nest2} combined with a scheme for restarting the algorithm has a convergence rate upper bound $f(y_k) - f(x^\ast)= O(\exp({-\mu^{\frac{1}{p}}\delta k}))$. We can consider this algorithm an accelerated method given the original descent method satisfies $f(y_k) - f(x^\ast)= O(\exp({-\mu^{\frac{1}{p-1}}\tilde \delta k}))$ under the same condition, where $\tilde \delta^{p-1} = \delta^{p}$. See \cref{app:restart} for details.
 \end{remark}
To summarize, it is sufficient to establish conditions under which an algorithm is a $\delta$-descent algorithm of order $p$ in order to (1) obtain a convergence rate and (2) accelerate the algorithm (in most cases). 

\paragraph{Accelerated rescaled gradient descent (Nesterov-style)}
Using~\cref{eq:nest1} we accelerate RGD. 
\begin{algorithm}[H]
\caption{Nesterov-style accelerated rescaled gradient descent.}
\label{eq:argd}
\begin{algorithmic}[1]
\Require{$f$ satisfies~\cref{eq:ss} and $h$ satisfies $D_h(x,y) \geq \frac{1}{p}\|x - y\|^p$}
	\State Set $x_0 = z_0$, $A_k = (\delta/p)^p k^{(p)}$, $\alpha_k =  \frac{A_{k+1} - A_k}{\delta}$, $\tau_k = \frac{\alpha_k}{A_{k+1}}$, and $\delta^{\frac{p}{p-1}} = \eta^{\frac{1}{p-1}}/2$.\\ {\bf for} $k = 1, \dots, K$ {\bf do}
\State $x_{k} =  \delta \tau_k z_k + (1 - \delta \tau_k) y_k$
\State $z_{k+1} = \arg \min_{z \in \mathcal{X}}\left\{ \alpha_k\langle \nabla f(x_k), z\rangle +  \frac{1}{\delta}D_h(z, z_k)\right\}$ 
\State $y_{k+1} =x_k -  \eta^{\frac{1}{p-1}}B^{-1}\nabla f(x_k)/\|\nabla f(x_k)\|_{\ast}^{\frac{p-2}{p-1}} $%
\State \Return $y_K$.
\end{algorithmic}
\end{algorithm}
We summarize the performance of \cref{eq:argd} in the following Corollary to Theorems~\ref{lem:1} and~\ref{prop:8}:
\begin{theorem}
Suppose $f$ is convex and strongly smooth of order $1<p<\infty$  with constants $0<L_1,\dots,L_p<\infty$. Also suppose $\eta$ satisfies~\cref{eq:step}. Then \cref{eq:argd} satisfies the convergence rate upper bound~\cref{eq:rate1}.
\end{theorem}

\subsection{Monteiro-Svaiter acceleration of descent algorithms}
Recently, \citet{Svaiter} have introduced an alternative framework for accelerating descent methods, which is similar to Nesterov's scheme but includes a line search step. This framework was further generalized by several more recent concurrent works~\cite{Gasnikov,Optimal,Bubeck} who demonstrate that higher-order tensor method~\cref{eq:hgd} with the addition of a line search step obtains a convergence rate upper bound $f(y_k) - f(x^\ast) = O( 1/k^{\frac{3p-2}{2}})$. When $p =2$, this rate matches that of the Nesterov-style acceleration framework, but for $p>2$ it is better. In this section, we present a novel, generalized version of the Monteiro-Svaiter accleration framework. In particular, {\em we use a simple Lyapunov analysis to generalize the framework} and show that many other descent methods of order $p$ can be accelerated in it, including the proximal method~\cref{eq:prox}, RGD~\cref{eq:rgd} and universal tensor methods.
\begin{theorem}
\label{prop:MS}
Suppose $h$ is satisfies the condition $B\preceq \nabla^2 h $. Consider sequence~\cref{eq:nest1} where in addition, we add a line search step which ensures the inequalities
\begin{subequations}\label{eq:ms1}
\begin{talign}
a \leq \frac{\lambda_{k+1}}{\delta^{\frac{3p-2}{2}}}\|y_{k+1} - x_k\|^{p-2}&\leq b, \quad 0<a<b \qtext{and}\label{eq:ms}\\
\|y_{k+1} - x_k +\lambda_{k+1} \nabla f(y_{k+1})\|&\leq \frac{1}{2}\| y_{k+1}  - x_k\| \label{eq:ms2}
\end{talign}
\end{subequations}
hold for the pair $(\lambda_{k+1}, y_{k+1})$, where $\lambda_{k+1} = \delta^2\alpha_k^2/A_{k+1}$ . Then the composite sequence satisfies:
\begin{talign}\label{eq:rate2} 
f(y_k) - f(x^\ast) \leq \frac{p^{\frac{3p-2}{2}}D_h(x^\ast, x_0)^{\frac{p}{2}}}{(\delta k)^{\frac{3p-2}{2}}} = O\left(1/(\delta k)^{\frac{3p-2}{2}}\right).
\end{talign}
\end{theorem}
The proof of \cref{prop:MS} is in \cref{App:MSAcc}.
All the aforementioned concurrent works have demonstrated that the higher-order gradient method ($\nu = 1$) with the addition line search step satisfies~\cref{eq:ms1}. We show the same is true of the proximal method~\cref{eq:prox}, rescaled gradient descent~\cref{eq:rgd} and universal higher-order tensor methods. See \cref{app:thmM-ARGD} for details. We conjecture that all methods that satisfy conditions~\cref{eq:ms} and~\cref{eq:ms2} are descent methods of order $p$ with an additional line search step.%
\begin{remark}[Restarting for improved accelerated linear rate] If, in addition, $f$ is $\mu$-gradient dominated of order $p$, then~\cref{eq:ms1} combined with a scheme for restarting the algorithm satisfies the convergence rate upper bound $f(y_k) - f(x^\ast) = O(\exp(-\mu^{\frac{2}{3p-2}}\delta k))$. See \cref{app:restart} for details.
\end{remark}

\subsection{Accelerating rescaled gradient descent (Monteiro-Svaiter-style)}
 Monteiro-Svaiter accelerated rescaled gradient descent is the following algorithm.

\begin{algorithm}[H]
\caption{Monteiro-Svaiter-style accelerated rescaled gradient descent.}
\label{alg:argd2}
\begin{algorithmic}[1]
\Require{$f$ is {\em strongly smooth of order $1<p<\infty$} and $h$ satisfies $B \preceq \nabla^2 h $.}
	\State Set $x_0 = z_0 =0$, $A_0 = 0$, $\delta^{\frac{3p-2}{2}} = \eta$,\,$\eta^{\frac{1}{p-1}} \leq \min\{ \frac{2}{5p}, 1/(2\sum_{m=2}^p \frac{L_m}{m!})\}$\\
	{\bf for} $k = 1, \dots, K$ {\bf do}
	\begin{subequations}
		\State Choose $\lambda_{k+1}$ (e.g. by line search) such that 
	$\frac{3}{4} \leq \frac{\lambda_{k+1} \|y_{k+1} -  x_k\|^{p-2}}{\eta} \leq \frac{5}{4}$,
	where 
	\begin{align*} 
	y_{k+1} = x_k - \eta^{\frac{1}{p-1}} \frac{\nabla f(x_k)}{\|\nabla f(x_k)\|_\ast^{\frac{p-2}{p-1}}},
	\end{align*}
	\end{subequations}
	and $\alpha_{k} =\frac{ \lambda_{k+1} + \sqrt{\lambda_{k+1} + 4A_k\lambda_{k+1}}}{2\delta}$, $A_{k+1} = \delta \alpha_{k}  + A_k$, $\tau_k = \frac{\alpha_k}{A_{k+1}}$ (so that $\lambda_{k+1} = \frac{\delta^2 \alpha_k^2}{A_{k+1}}$) and 
		$$\textstyle x_k = \delta \tau_k  z_k + (1-\delta \tau_k) y_k.$$

	\State Update $z_{k+1} = \arg \min_{z\in \mathcal{X}}\left\{ \alpha_k \langle \nabla f(y_{k+1}), z\rangle +   \frac{1}{\delta}D_h(z, z_k)\right\}$
\State \Return $y_K$.
\end{algorithmic}
\end{algorithm}
We summarize results on performance of \cref{alg:argd2} in the following corollary to \cref{prop:MS}:
\begin{theorem}\label{thm:M-ARGD} Assume $f$ is convex and strongly smooth of order $1 < p <\infty$ with constants $0<L_1, \dots, L_p<\infty$. Then \cref{alg:argd2} satisfies the convergence rate upper bound~\cref{eq:rate2}.
\end{theorem}

\section{Related Work}
Our acceleration framework is similar in spirit to a number of acceleration frameworks in the literature (e.g., \citet{coupling}, \citet{Lessard14}, \citet{Zaid}, \citet{jelena}) but applies more generally to descent methods of order $p>2$. In particular, the present framework builds off of the framework proposed by \citet{Wilson}, but it (1) makes the connection to descent methods more explicit and (2) incorporates a generalization and Lyapunov analysis of the Monteiro-Svaiter acceleration framework. These manifold generalizations crucially allow us to propose RGD and accelerated RGD, which has superior theoretical and empirical performance to several existing methods on strongly smooth functions. %
\section{Examples and Numerical Experiments}
\label{sec:exam}
We compare our result to several recent works that have shown that for some function classes, more intuitive first-order algorithms outperform gradient descent. In particular, both \citet{Runge} and \citet{Hamiltonian} obtain first-order algorithms by applying integration techniques to second-order ODEs. When the objective function is sufficiently smooth, both show their algorithm outperforms (accelerated) gradient descent.
We show that %
Algorithms~\ref{eq:argd} and~\ref{alg:argd2} achieves fast performance in theory and in practice on similar objectives. %
\paragraph{Runge-Kutta}
\citet{Runge} show that if one applies an $s$-th order Runge-Kutta integrator to a family of second-order dynamics, then the resulting algorithm\footnote{which requires at least $s$ gradient evaluations per iteration} achieves a convergence rate $f(x_k) - f(x^\ast) =O(1/k^{\frac{ps}{s-1}})$\footnote{this matches the rate of \cref{eq:argd} in the limit  $s \rightarrow \infty$, where $s$ is the order of the integrator.} provided the function meet the following two conditions: (1)
$f$ satisfies the {\em gradient lower bound} of order $p \ge 2$, which means for all $m = 1,\dots,p-1$,
\begin{talign}\label{eq:rung1}
\textstyle f(x)-f(x^\ast) \ge \frac{1}{C_m} \|\nabla^m f(x)\|^{\frac{p}{p-m}} ~~~ \forall \, x \in \R^n
\end{talign}
for some constants $0 < C_1, \dots, C_{p-1} < \infty$; and (2) for $s \geq p$ and $M>0$, $f$ is $(s+2)$-times differentiable and $\|\nabla^{(i)}f(x)\|\leq M$ for $i=p, p+1,\dots, s+2$. 
One can show that if $f$ is strongly smooth of order $p$, then $f$ satisfies the gradient lower bound of order $p$. The details of this result is in~\ref{supp:exp}. While are unable to prove that condition~\cref{eq:rung1} is equivalent to strong smoothness, we have yet to find an example of a function that satisfies~\cref{eq:rung1} and is not strongly smooth.%
\paragraph{Hamiltonian Descent} \citet{Hamiltonian} show explicit integration techniques applied to conformal Hamiltonian dynamics converge at a fast linear rate for a function class larger than gradient descent. The method entails finding a kinetic energy map that upper bounds the dual of the function. All examples for which we can compute such a map given by~\cite{Hamiltonian} are uniformly convex and gradient dominated functions; therefore, simply rescaling the gradient for these examples ensures a linear rate. 

\begin{figure}[htb] %
\centering
\begin{subfigure}{0.5\textwidth}
\includegraphics[width=\linewidth]{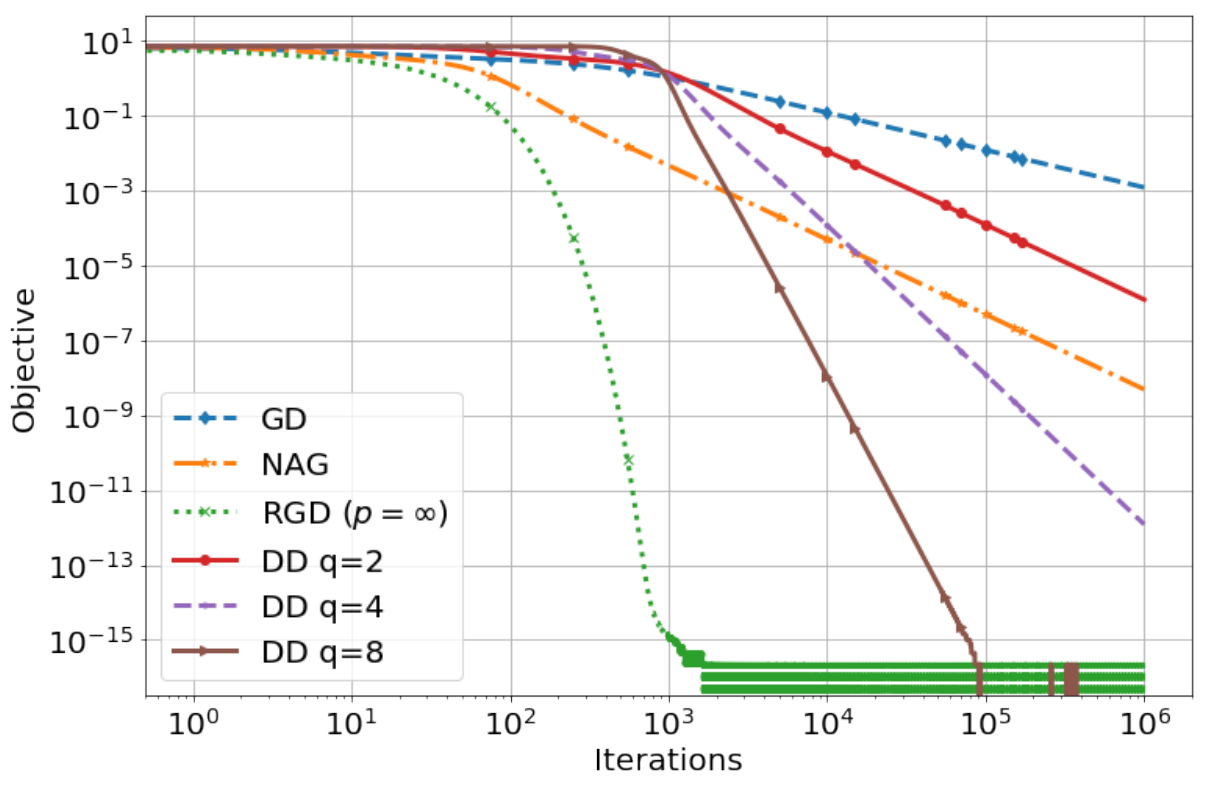}
\caption{Example~\ref{exam:1}: Logistic loss (Iteration)} \label{fig:a}
\end{subfigure}%
\begin{subfigure}{0.5\textwidth}
\includegraphics[width=\linewidth]{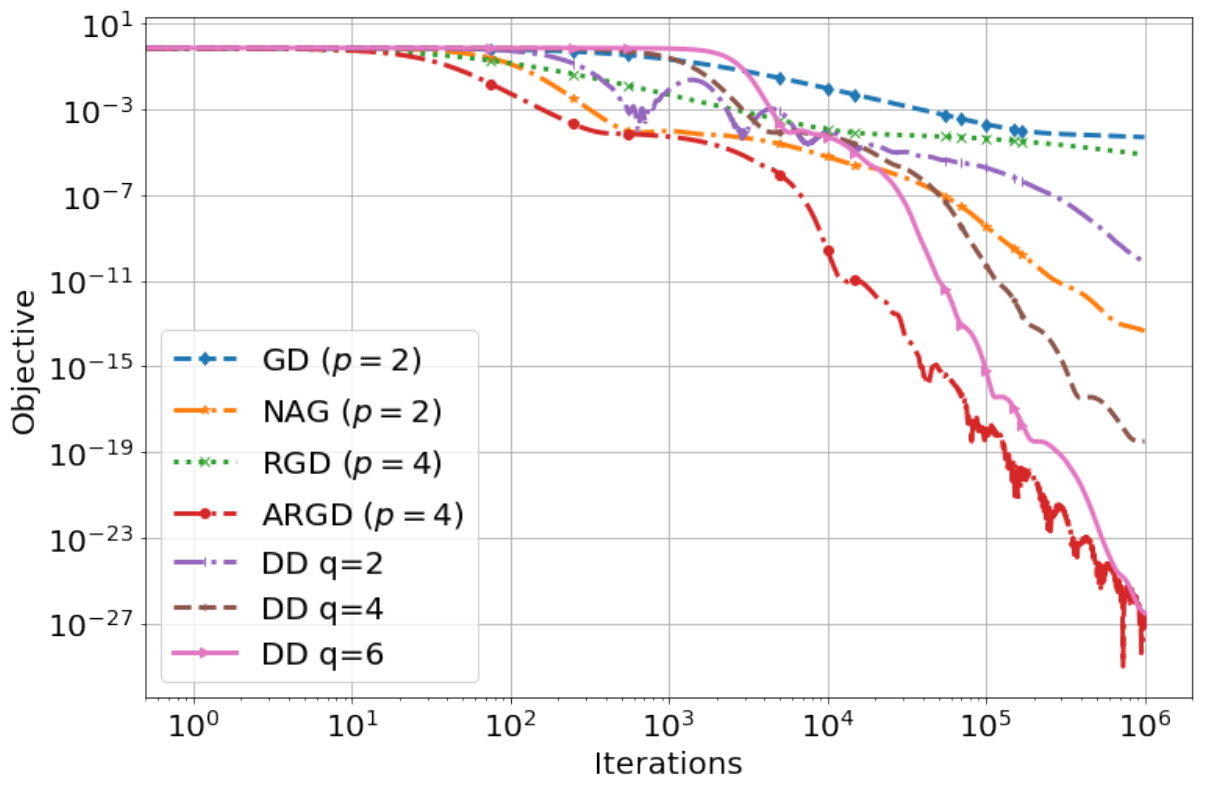}
\caption{Example~\ref{exam:2}$:\ell_4$ loss (Iteration)} \label{fig:b}
\end{subfigure}
\begin{subfigure}{0.5\textwidth}
\includegraphics[width=\linewidth]{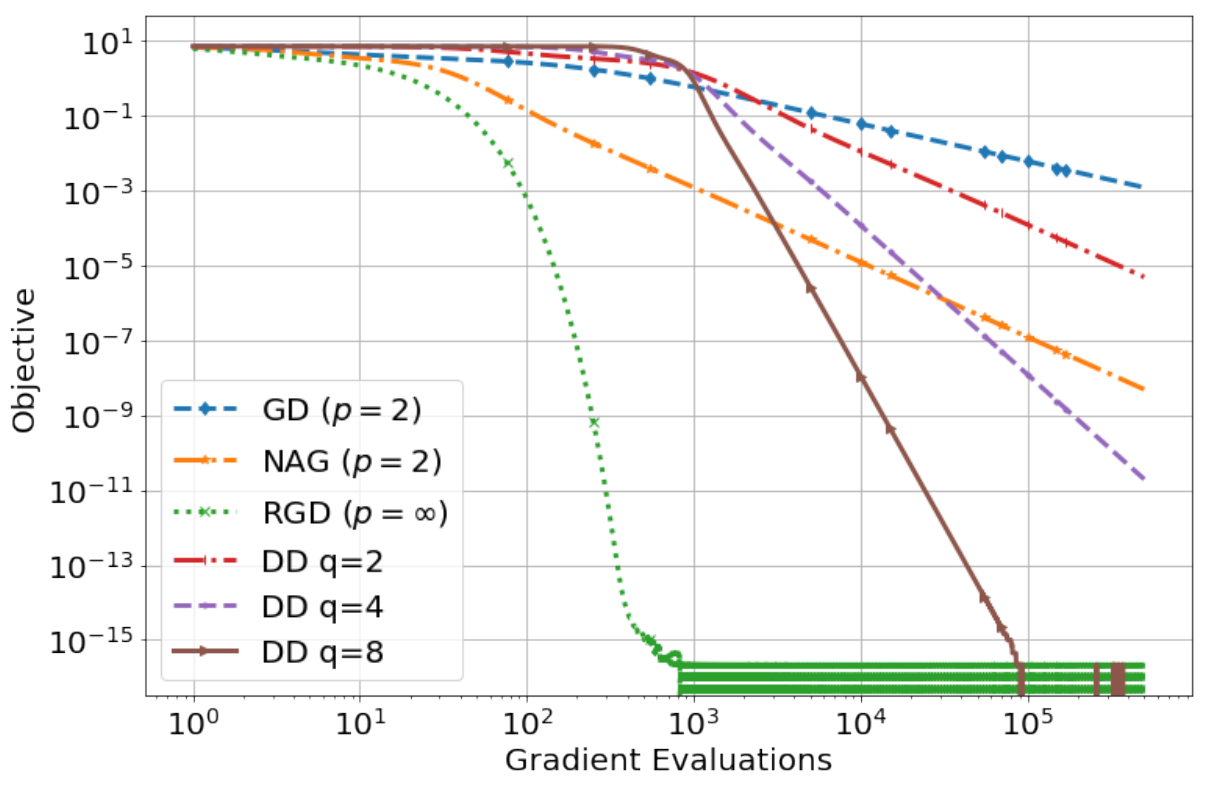}
\caption{Example~\ref{exam:1}: Logistic loss (Gradient)} \label{fig:c}
\end{subfigure}%
\begin{subfigure}{0.5\textwidth}
\includegraphics[width=\linewidth]{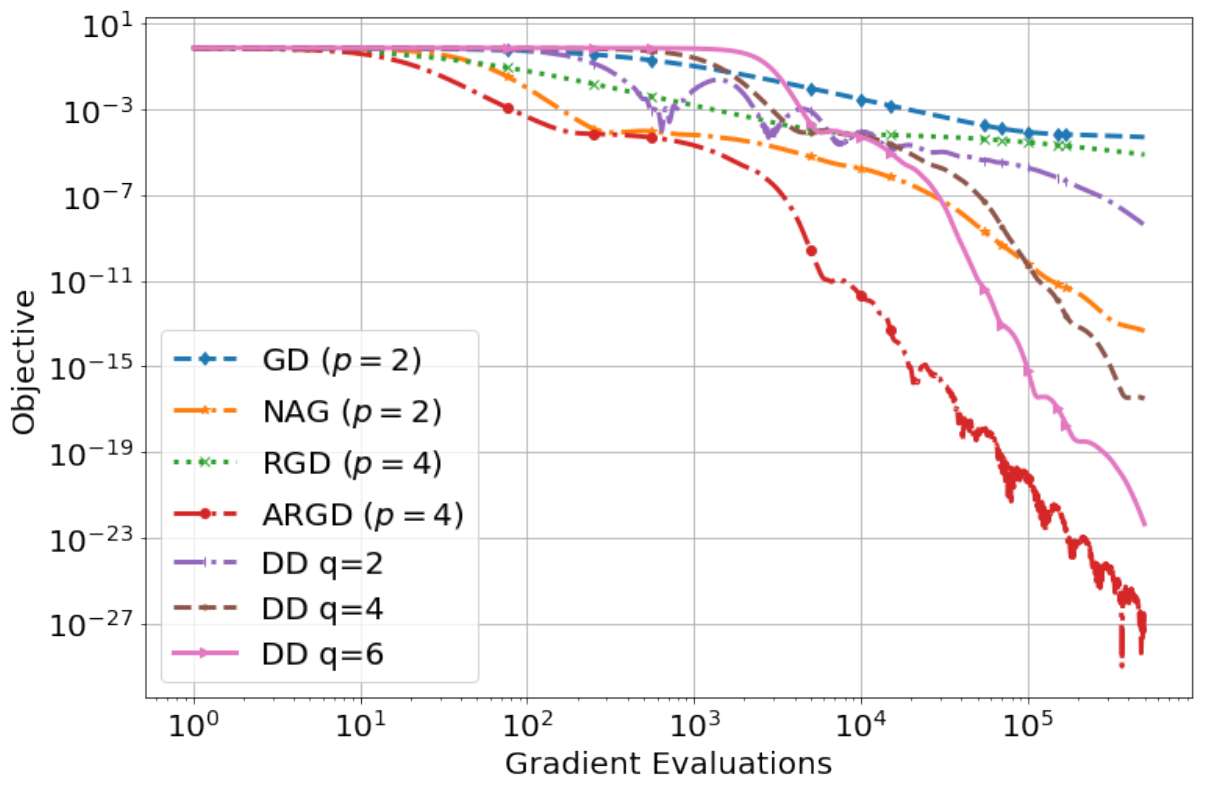}
\caption{Example~\ref{exam:2}: $\ell_4$ loss (Gradient)} \label{fig:d}
\end{subfigure}
\medskip
\begin{center}
\begin{subfigure}{0.5\textwidth}
\includegraphics[width=\linewidth]{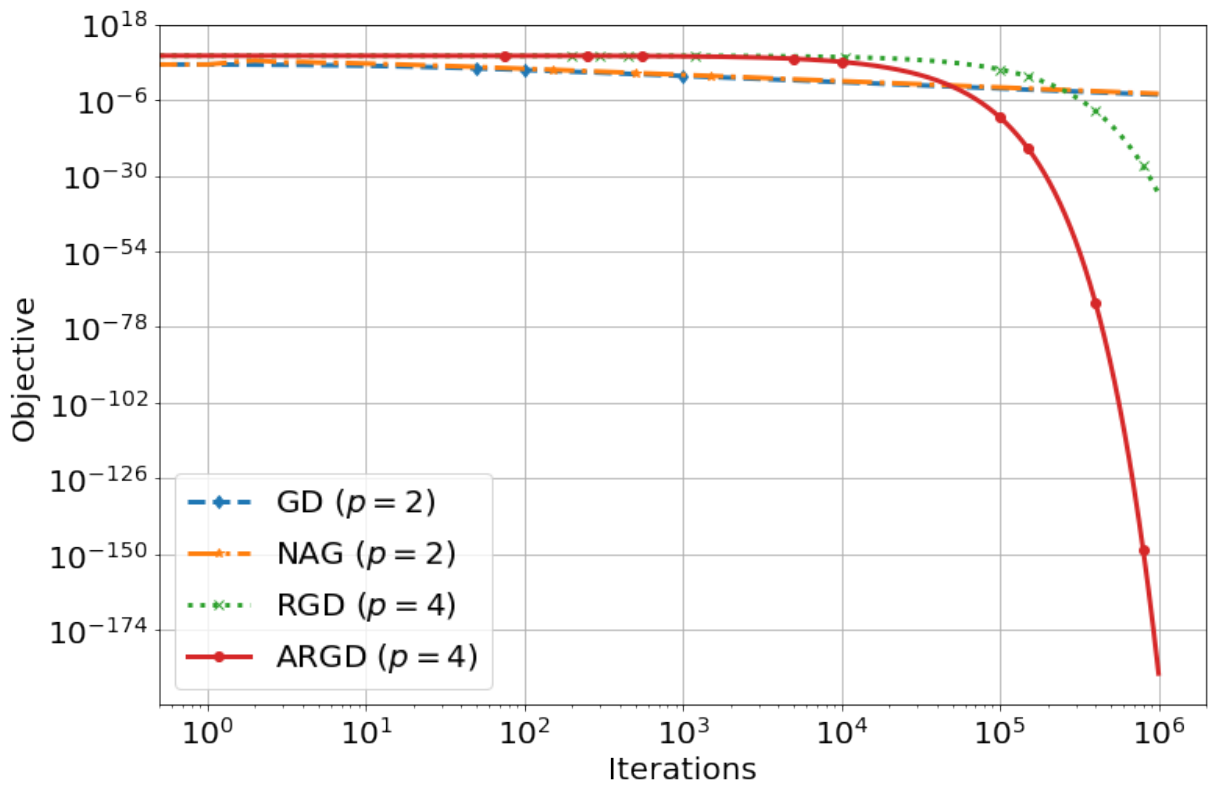}
\caption{Example~\ref{exam:5}:  Hamiltonian function} \label{fig:e}
\end{subfigure}%
\end{center}
    \caption{Experimental results comparing RGD and accelerated RGD (ARGD) to gradient descent (GD), Nesterov accelerated GD (NAG) and Runge-Kutta (DD).  The plots for Runge-Kutta use an $s = 2$ integrator which requires two gradient evaluations per iteration.  
 Where relevant, we plot both iterations (\cref{fig:a,fig:b}) and gradient evaluations (\cref{fig:c,fig:d}).}
\end{figure}

\subsection{Examples}
We provide several examples of strongly smooth functions in machine learning (see \cref{app:example} for details).
\begin{example} \label{exam:2}The $\ell_p$ loss function 
\begin{talign}%
f(x) = \frac{1}{p}\|A x - b\|_p^p,
\end{talign}
shown by \citet{Runge} to satisfy~\cref{eq:rung1} of order $p$, is strongly smooth of order $p$. \end{example}
\begin{example}\label{exam:1}
The logistic loss
\begin{talign}
f(x) =  \log (1 + e^{-yw^\top x}),
\end{talign}
shown by \citet{Runge} to satisfy~\cref{eq:rung1} of order $p = \infty$, is strongly smooth of order $p=\infty$. 
\end{example}
\begin{example}\label{exam:3}
The GLM loss, 
\begin{talign}
f(x) = \frac{1}{2} (y - \phi(x^\top w))^2 \qtext{for} \phi(r) = 1/(1+e^{-r}),\quad y \in \{0,1\}, \qtext{ and } w \in \R^d,
\end{talign} 
studied by \citet{Hazan} is strongly smooth of order $p=3$.
\end{example}
\begin{example}\label{exam:4}
The $\ell_2$ loss to the $p$-th power
\begin{talign}
f(x) = \frac{1}{p}\|A x - b\|_2^p,
\end{talign}
for which Hamiltonian descent~\cite{Hamiltonian} obtains a linear rate, is strongly smooth and gradient dominated of order $p$.
\end{example}
\begin{example}\label{exam:5} The loss function,
\begin{talign}
f(x) = (x^{(1)} + x^{(2)})^4 + \frac{1}{16}(x^{(1)} - x^{(2)})^4, 
\end{talign}
for which Hamiltonian descent~\cite{Hamiltonian} obtains a linear rate, is strongly smooth and gradient dominated of order $p=4$.
\end{example}
\subsection{Experiments}
In this section, we perform a series of numerical experiments to compare the performance of ARGD (\cref{eq:argd}) with gradient descent (GD), Nesterov accelerated GD (NAG), and the state-of-the-art Runge-Kutta algorithms of~\citet{Runge} (DD) on the logistic loss $f(x) = \sum_{i=1}^{10} \log(1 + e^{-w_i^\top x y_i})$, the $\ell_4$ loss $f(x) = \frac{1}{4}\|A x - b\|_4^4$, and the Hamiltonian descent loss (\cref{exam:5}).  For the logistic and $\ell_4$ losses, we use the same code, plots, and experimental methodology of \citet{Runge} (including data and step-size choice), adding to it (A)RGD. Specifically, for~\cref{fig:a}-\cref{fig:d}, the entries of $W \in \R^{10 \times 10} $ and $A \in \R^{10 \times 10}$ are i.i.d. standard Gaussian, and the first five entries of $y$ (and $b$) are valued 0 while the rest are 1.
\cref{fig:e} shows the performance of A(RGD), GD, and NAG on the Hamiltonian objective studied by~\cite{Hamiltonian}; for~\cref{fig:e}, the largest step-size was chosen subject to the algorithm not diverging.
 For each experiment, a simple implementation of (A)RGD significantly outperforms the Runge-Kutta algorithm (DD), GD and NAG. The code for these experiments can be found here: \url{https://github.com/aswilson07/ARGD.git}.
\section{Additional Results and Discussion}
\label{sec:4}
This paper establishes broad conditions under which an algorithm will converge and its performance can be accelerated by adding momentum. We use these conditions to introduce (accelerated) rescaled gradient descent for strongly smooth functions, and showed it outperforms several recent first-order methods that have been introduced for optimizing smooth functions in machine learning. 

There are (at least) two simple extensions of our framework. First, %
 an analogous framework can be established for {\bf (accelerated) $\delta$-coordinate descent methods} of order $p$.  As an application, we introduce (accelerated) rescaled coordinate descent for functions that are strongly smooth along each coordinate direction of the gradient. We provide details in \cref{app:Cor}. Second, 
with our generalization of the Monteiro-Svaiter framework, we derive {\bf optimal univeral tensor methods} for functions whose $(p-1)$-st gradients are $\nu$-H\"older-smooth which achieve the upper bound $f(y_k) - f(x^\ast) = O(1/k^{\frac{3\tilde{p}-2}{2}})$ where $\tilde p = p -1 + \nu$. The matching lower bound for this class of functions was recently established by~\cite{Nesterov19}. We present this result in \cref{app:Uni}.

There are several possible directions for future work. We know that certain simple operations preserve convexity (e.g., addition), but what operations preserve strong smoothness?  Understanding this could allow us to construct more complex examples of strongly smooth functions. %
 Our results reveal an interesting hierarchy of smoothness assumptions which lead to methods that converge quickly; exploring this more is of significant interest. Finally, extending our analysis to the stochastic or manifold setting, studying the use of  variance reduction techniques, and introducing other $\delta$-decent algorithms of order $p$ are all interesting directions for future work. 

\subsubsection*{Acknowledgments}
We would like to thank Jingzhao Zhang for providing us access to his code. 
\bibliography{refs}

\begin{thebibliography}{36}
\providecommand{\natexlab}[1]{#1}
\providecommand{\url}[1]{\texttt{#1}}
\expandafter\ifx\csname urlstyle\endcsname\relax
  \providecommand{\doi}[1]{doi: #1}\else
  \providecommand{\doi}{doi: \begingroup \urlstyle{rm}\Url}\fi

\bibitem[{Allen Zhu} and Orecchia(2017)]{coupling}
Zeyuan {Allen Zhu} and Lorenzo Orecchia.
\newblock Linear coupling: An ultimate unification of gradient and mirror
  descent.
\newblock In \emph{8th Innovations in Theoretical Computer Science Conference,
  {ITCS} 2017, January 9-11, 2017, Berkeley, CA, {USA}}, pages 3:1--3:22, 2017.

\bibitem[Amari(1998)]{Amari1988}
Shun-Ichi Amari.
\newblock Natural gradient works efficiently in learning.
\newblock \emph{Neural Computation}, pages 251--276, 1998.

\bibitem[Baes(2009)]{Baes09}
Michel Baes.
\newblock Estimate sequence methods: Extensions and approximations, August
  2009.

\bibitem[Betancourt et~al.(2018)Betancourt, Jordan, and
  Wilson]{BetancourtJordanWilson}
Michael Betancourt, Michael Jordan, and Ashia Wilson.
\newblock On symplectic optimization.
\newblock Arxiv preprint arXiv{1802.03653}, 2018.

\bibitem[Bubeck et~al.()Bubeck, Jiang, Lee, Li, and Sidford]{Bubeck}
S{\'e}bastien Bubeck, Qijia Jiang, Yin~Tat Lee, Yuanzhi Li, and Aaron Sidford.
\newblock Near-optimal method for highly smooth convex optimization.
\newblock In \emph{Proceedings of the Thirty-Second Conference on Learning
  Theory}, volume~99 of \emph{Proceedings of Machine Learning Research}, pages
  492--507, Phoenix, USA, 25--28 Jun . PMLR.

\bibitem[Carmon et~al.(2017)Carmon, Duchi, Hinder, and Sidford]{Yair}
Yair Carmon, John~C Duchi, Oliver Hinder, and Aaron Sidford.
\newblock Lower bounds for finding stationary points ii: First-order methods.
\newblock Arxiv preprint arXiv:1711.0084, 2017.

\bibitem[Chen and Teboulle(1993)]{ChenTeboulle93}
Gong Chen and Marc Teboulle.
\newblock Convergence analysis of a proximal-like minimization algorithm using
  {B}regman functions.
\newblock \emph{SIAM Journal of Optimization}, 3\penalty0 (3):\penalty0
  538--543, 1993.

\bibitem[Diakonikolas and Orecchia(2018)]{jelena}
Jelena Diakonikolas and Lorenzo Orecchia.
\newblock Accelerated extra-gradient descent: {A} novel accelerated first-order
  method.
\newblock In \emph{9th Innovations in Theoretical Computer Science Conference,
  {ITCS} 2018, January 11-14, 2018, Cambridge, MA, {USA}}, pages 23:1--23:19,
  2018.

\bibitem[Gasnikov et~al.(2019)Gasnikov, Dvurechensky, Gorbunov, Vorontsova,
  Selikhanovych, and Uribe]{Gasnikov}
Alexander Gasnikov, Pavel Dvurechensky, Eduard Gorbunov, Evgeniya Vorontsova,
  Daniil Selikhanovych, and C\'esar~A. Uribe.
\newblock Optimal tensor methods in smooth convex and uniformly convex
  optimization.
\newblock In \emph{Proceedings of the Thirty-Second Conference on Learning
  Theory}, pages 1374--1391, Phoenix, USA, 25--28 Jun 2019. PMLR.

\bibitem[Grapiglia and Nesterov(2019)]{Nesterov19}
G.N Grapiglia and Yu. Nesterov.
\newblock Tensor methods for minimizing functions with h\"older continuous
  higher-order derivatives.
\newblock Arxiv preprint arXiv{1904.12559}, April 2019.

\bibitem[Hazan et~al.(2015)Hazan, Levy, and Shalev{-}Shwartz]{Hazan}
Elad Hazan, Kfir~Y. Levy, and Shai Shalev{-}Shwartz.
\newblock Beyond convexity: Stochastic quasi-convex optimization.
\newblock In \emph{Advances in Neural Information Processing Systems 28: Annual
  Conference on Neural Information Processing Systems 2015, December 7-12,
  2015, Montreal, Quebec, Canada}, pages 1594--1602, 2015.

\bibitem[Jiang et~al.(2018)Jiang, Wang, and Zhang]{Optimal}
B.~Jiang, H.~Wang, and S.~Zhang.
\newblock An optimal high-order tensor method for convex optimization.
\newblock \emph{Arxiv preprint arXiv:1812.06557}, 2018.

\bibitem[Krichene et~al.(2015)Krichene, Bayen, and Bartlett]{KricheneStoch}
Walid Krichene, Alexandre Bayen, and Peter~L Bartlett.
\newblock Accelerated mirror descent in continuous and discrete time.
\newblock In C.~Cortes, N.~D. Lawrence, D.~D. Lee, M.~Sugiyama, and R.~Garnett,
  editors, \emph{Advances in Neural Information Processing Systems 28}, pages
  2845--2853. Curran Associates, Inc., 2015.

\bibitem[Lessard et~al.(2016)Lessard, Recht, and Packard]{Lessard14}
Laurent Lessard, Benjamin Recht, and Andrew Packard.
\newblock Analysis and design of optimization algorithms via integral quadratic
  constraints.
\newblock \emph{SIAM Journal on Optimization}, 26\penalty0 (1):\penalty0
  57--95, 2016.

\bibitem[Lin et~al.(2017)Lin, Mairal, and Harchaoui]{Zaid}
Hongzhou Lin, Julien Mairal, and Za{\"{\i}}d Harchaoui.
\newblock Catalyst acceleration for first-order convex optimization: from
  theory to practice.
\newblock \emph{Journal of Machine Learning Research}, 18:\penalty0
  212:1--212:54, 2017.

\bibitem[\L{}ojasiewicz(1963)]{Loj64}
S.~\L{}ojasiewicz.
\newblock A topological property of real analytic subsets (in french).
\newblock In \emph{Coll. du CNRS, Les \'equations aux d\'eri\'vees partielles},
  pages 87-- 89, 1963.

\bibitem[Maddison et~al.(2018)Maddison, Paulin, Teh, O'Donoghue, and
  Doucet]{Hamiltonian}
Chris~J. Maddison, Daniel Paulin, Yee~Whye Teh, Brendan O'Donoghue, and Arnaud
  Doucet.
\newblock Hamiltonian descent methods.
\newblock Arxiv preprint arXiv{1809.05042}, 2018.

\bibitem[Monteiro and Svaiter(2013)]{Svaiter}
Renato D.~C. Monteiro and Benar~Fux Svaiter.
\newblock An accelerated hybrid proximal extragradient method for convex
  optimization and its implications to second-order methods.
\newblock \emph{SIAM Journal on Optimization}, 23\penalty0 (2):\penalty0
  1092--1125, 2013.

\bibitem[Moreau(1965)]{prox}
Jean~Jacques Moreau.
\newblock Proximit\'e et dualit\'e dans un espace {Hilbertien}.
\newblock \emph{Bulletin de la Soci\'et\'e Math\'ematique de France},
  93:\penalty0 273--299, 1965.

\bibitem[Nemirovskii and Yudin(1983)]{NemirovskiiYudin}
Arkadi Nemirovskii and David Yudin.
\newblock \emph{Problem {C}omplexity and {M}ethod {E}fficiency in
  {O}ptimization}.
\newblock John Wiley \& Sons, 1983.

\bibitem[Nesterov(2018)]{Nesterov18}
Y.~Nesterov.
\newblock Implementable tensor methods in unconstrained convex optimization.
\newblock Core discussion papers, 2018.
\newblock URL \url{https://ideas.repec.org/p/cor/louvco/2018005.html}.

\bibitem[Nesterov(1983)]{Nesterov83}
Yurii Nesterov.
\newblock A method of solving a convex programming problem with convergence
  rate ${O}(1/k^2)$.
\newblock \emph{Soviet Mathematics Doklady}, 27\penalty0 (2):\penalty0
  372--376, 1983.

\bibitem[Nesterov(2004)]{Nesterov04}
Yurii Nesterov.
\newblock \emph{Introductory Lectures on Convex Optimization: A Basic Course}.
\newblock Applied Optimization. Kluwer, Boston, 2004.

\bibitem[Nesterov(2005)]{Nesterov05}
Yurii Nesterov.
\newblock Smooth minimization of non-smooth functions.
\newblock \emph{Mathematical Programming}, 103\penalty0 (1):\penalty0 127--152,
  2005.

\bibitem[Nesterov(2008)]{Nesterov08}
Yurii Nesterov.
\newblock Accelerating the cubic regularization of {N}ewton's method on convex
  problems.
\newblock \emph{Mathematical Programming}, 112\penalty0 (1):\penalty0 159--181,
  2008.
\newblock ISSN 0025-5610.

\bibitem[Nesterov and Polyak(2006)]{Nesterov06}
Yurii Nesterov and Boris~T. Polyak.
\newblock Cubic regularization of {N}ewton's method and its global performance.
\newblock \emph{Mathematical Programming}, 108\penalty0 (1):\penalty0 177--205,
  2006.

\bibitem[Polyak(1964)]{Polyak1964}
Boris~T. Polyak.
\newblock Some methods of speeding up the convergence of iteration methods.
\newblock \emph{{USSR} Computational Mathematics and Mathematical Physics},
  4\penalty0 (5):\penalty0 1--17, 1964.

\bibitem[Schropp and Singer(2000)]{Schropp}
J.~Schropp and I.~Singer.
\newblock A dynamical systems approach to constrained minimization.
\newblock \emph{Numerical Functional Analysis and Optimization}, 21\penalty0
  (3-4):\penalty0 537--551, 2000.

\bibitem[Shi et~al.(2018)Shi, Du, Jordan, and Su]{SuJordan}
Bin Shi, Simon Du, Michael Jordan, and Weiji Su.
\newblock Understanding the acceleration phenomenon via high-resolution
  differential equations.
\newblock Arxiv preprint arXiv{1810.08907}, November 2018.

\bibitem[Su et~al.(2014)Su, Boyd, and Cand\`{e}s]{SuBoydCandes14}
Weijie Su, Stephen Boyd, and Emmanuel~J. Cand\`{e}s.
\newblock A differential equation for modeling {N}esterov's accelerated
  gradient method: Theory and insights.
\newblock In \emph{Advances in Neural Information Processing Systems (NIPS)
  27}, 2014.

\bibitem[Sundaramoorthi and Yezzi(2018)]{PDE}
Ganesh Sundaramoorthi and Anthony~J. Yezzi.
\newblock Variational {PDEs} for acceleration on manifolds and application to
  diffeomorphisms.
\newblock In \emph{Advances in Neural Information Processing Systems 31: Annual
  Conference on Neural Information Processing Systems 2018, NeurIPS 2018, 3-8
  December 2018, Montr{\'{e}}al, Canada.}, pages 3797--3807, 2018.

\bibitem[Wibisono(2018)]{Andre}
Andre Wibisono.
\newblock Sampling as optimization in the space of measures: The {Langevin}
  dynamics as a composite optimization problem.
\newblock In \emph{Conference On Learning Theory, {COLT} 2018, Stockholm,
  Sweden, 6-9 July 2018.}, pages 2093--3027, 2018.

\bibitem[Wibisono and Wilson(2015)]{Accel2}
Andre Wibisono and Ashia Wilson.
\newblock On accelerated methods in optimization.
\newblock Arxiv preprint arXiv{1509.03616}, 2015.

\bibitem[Wibisono et~al.(2016)Wibisono, Wilson, and Jordan]{Acceleration}
Andre Wibisono, Ashia~C. Wilson, and Michael~I. Jordan.
\newblock A variational perspective on accelerated methods in optimization.
\newblock \emph{Proceedings of the National Academy of Sciences}, 113\penalty0
  (47):\penalty0 E7351--E7358, 2016.

\bibitem[Wilson et~al.(2016)Wilson, Recht, and Jordan]{Wilson}
Ashia Wilson, Benjamin Recht, and Michael Jordan.
\newblock A {Lyapunov }analysis of momentum methods in optimization.
\newblock Arxiv preprint arXiv{1611.02635}, November 2016.

\bibitem[Zhang et~al.(2018)Zhang, Mokhtari, Sra, and Jadbabaie]{Runge}
Jingzhao Zhang, Aryan Mokhtari, Suvrit Sra, and Ali Jadbabaie.
\newblock Direct {Runge-Kutta} discretization achieves acceleration.
\newblock In S.~Bengio, H.~Wallach, H.~Larochelle, K.~Grauman, N.~Cesa-Bianchi,
  and R.~Garnett, editors, \emph{Advances in Neural Information Processing
  Systems 31}, pages 3904--3913. Curran Associates, Inc., 2018.

\end{thebibliography}
\bibliographystyle{plainnat}
\bibpunct{(}{)}{;}{a}{,}{,}

\newpage
\appendix
\allowdisplaybreaks

{\centering
\vbox{\hsize\textwidth
\linewidth\hsize \vskip 0.1in %
{\LARGE \textsc{Supplementary material to} \\
\bf Accelerating Rescaled Gradient Descent: \\Fast Minimization of Smooth Functions \par} %
\vspace{16pt}
\begin{tabular}[t]{c}\bf
Ashia C. Wilson\quad Lester Mackey \quad Andre Wibisono 
\end{tabular}\hfil\linebreak[0]\hfil
\vskip 0.3in}}
\section{Descent Flows}
The derivation and analysis of descent algorithms is inspired by {\em descent flows}. In this section we introduce and analyzed these family of dynamics.
\begin{definition} A dynamics is a {\em\bf descent flow of order $p$} if is satisfies:
\begin{align}\label{eq:ass1}
\frac{d}{dt} f(X_t) \leq -\|\nabla f(X_t)\|_{\ast}^{\frac{p}{p-1}},
\end{align}
for some $1<p\leq\infty$ and for all $0 \leq t \leq \infty$. 
\end{definition}
For dynamics that satisfy~\eqref{eq:ass1}, we obtain non-asymptotic convergence guarantees for non-convex, convex and gradient-dominated functions. 
We summarize our main results for descent curves of order $p$ in the following three theorems.
\begin{theorem}\label{thm:1}
 Suppose a dynamical system satisfies~\eqref{eq:ass1} for some $1<p\leq\infty $ and $f$ is differentiable.  Then the system satisfies:  
 \begin{align}\label{eq:cont1} 
\min_{0\leq s\leq t}\|\nabla f(X_s)\|_{\ast} = O \left(1/t^{\frac{p-1}{p}}\right).
\end{align}
\end{theorem}
\begin{theorem}\label{thm:2}
 Suppose a dynamical system satisfies~\eqref{eq:ass1} for some $1<p\leq\infty$ and $f$ is differentiable and convex with $R = \sup_{x:f(x)\leq f(x_0)} \|x - x^\ast \| < \infty$.  
 Then the system satisfies:
\begin{align}\label{eq:cont2}
f(X_t) - f(x^\ast) =
\begin{cases}
O\left(1/\Big(1+\frac{1}{Rp}{t^{\frac{p-1}{p}}}\Big)^{p}\right)&\text{ if }\, p < \infty  \\
O\left(e^{-t/R}\right)&\text{ if }\, p = \infty
 \end{cases}.
\end{align}
\end{theorem}
\begin{theorem}\label{thm:3}
 Suppose a dynamical system satisfies~\eqref{eq:ass1} for some $1<p\leq\infty$ and $f$ is differentiable and $\mu$-gradient dominated of order $p$.  Then the system satisfies: 
\begin{align}\label{eq:cont3}
f(X_t) - f(x^\ast) = O\left( e^{-\frac{p}{p-1}\mu^{\frac{1}{p-1}}  t}\right).
\end{align}
\end{theorem}
The proof of these results follows the same structure as the descent algorithms, with both relying on simple energy arguments. 
\subsection{Proofs}
To show~\eqref{eq:cont1},
we begin with the energy function $\mathcal{E}_t = f(X_t) - f(x^\ast)$. A quick calculation show:
\begin{align*}
\frac{d}{dt} \mathcal{E}_t = \frac{d}{dt} f(X_t) \overset{\eqref{eq:ass1}}{\leq} -\|\nabla f(X_t)\|_{\ast}^{\frac{p}{p-1}} 
\end{align*}
Integrating and rearranging gives the bound
\begin{align*}
 t \min_{0\leq s\leq t} \|\nabla f(X_s)\|^{\frac{p}{p-1}}  \leq \int_0^t  -\|\nabla f(X_t)\|_{\ast}^{\frac{p}{p-1}} dt \leq \E_0 - \E_t.
\end{align*}%
from which we can conclude~\eqref{eq:cont1}. To establish~\eqref{eq:cont2}, consider the energy function $\mathcal{E}_t = t^p(f(X_t) - f(x^\ast))$. We compute 
\begin{align*}
\frac{d}{dt} \mathcal{E}_t &= pt^{p-1}(f(X_t) - f(x^\ast))+ t^{p}\frac{d}{dt} f(X_t)\\
& \leq pt^{p-1}\langle \nabla f(X_t), x^\ast - X_t\rangle + t^{p}\frac{d}{dt} f(X_t)\\
&\overset{\eqref{eq:ass1}}{\leq}pt^{p-1}\langle \nabla f(X_t), x^\ast - X_t\rangle - t^{p}\|\nabla f(X_t)\|^{\frac{p}{p-1}}\\
&\leq \frac{1}{p}\|p(X_t - x^\ast)\|^p \leq p^{p-1} R^p.
\end{align*}
The first inequality uses the convexity of $f$ and the second inequality~\eqref{eq:ass1}. The third inequality uses the Fenchel-Young inequality $-\|s\|^{\frac{p}{p-1}} + \langle s,u \rangle \le - \frac{p-1}{p} \|s\|^{\frac{p}{p-1}} + \langle s,u \rangle \le \frac{1}{p}\|u\|^p$
 with $s = t^{p-1}\nabla f(X_t)$ and $ u =  p(x^\ast- X_t)$. The last step uses the fact that $\|X_t - x^\ast \| \leq R$ since \eqref{eq:ass1} implies the dynamical system is a descent method. Integrating allows us to obtain the statement $\E_t - \E_0 \leq tp^{p-1} R^p$, and subsequently, the upper bound
\begin{align*}
f(X_t) - f(x^\ast) \leq \frac{p^{p-1}R^p}{t^{p-1}},
\end{align*}
as desired. The last bound~\eqref{eq:cont3} uses the energy function $\mathcal{E}_t = f(X_t) - f(x^\ast)$ to establish
\begin{align*}
\frac{d}{dt} \mathcal{E}_t = \frac{d}{dt} f(X_t) \overset{\eqref{eq:ass1}}{\leq} -\|\nabla f(X_t)\|_{\ast}^{\frac{p}{p-1}} \leq \frac{p}{p-1}\mu^{\frac{1}{p-1}} \mathcal{E}_t.
\end{align*}
where the last inequality follows from the gradient dominated condition. We use the intuition from the bounds established for descent dynamics to derive analogous results for descent algorithms.
\section{Descent Algorithms}
We present proofs of results Section~\ref{sec:2}.
\label{App:A}
\subsection{Proof of Theorems~\ref{thm:11}-\ref{thm:33}}
We begin with detailed proofs of Theorems~\ref{thm:11}-\ref{thm:33}.
\subsubsection{Proof of Theorem~\ref{thm:11}}
By rearranging and summing~\eqref{eq:prog1}, we obtain
\begin{talign*} 
 \delta k \min_{j-k\leq s \leq j} \|\nabla f(x_{s})\|_{\ast}^{\frac{p}{p-1}} \leq \sum_{s={j-k}}^j \|\nabla f(x_{s})\|_{\ast}^{\frac{p}{p-1}} \delta\leq f(x_0) - f(x_k) \leq f(x_0)
\end{talign*}
where $j = k$ if the bound~\eqref{eq:p1} holds and $j = k+1$ if the bound~\eqref{eq:p2} holds. Rearranging the inequality yields the result in Theorem~\ref{thm:11}.

\subsubsection{Proof of Theorem~\ref{thm:22}}
Fix any $a > 0$, and define the positive increasing function $w_a(t) = (1+t/(ap))^p$, which satisfies $\frac{d}{dt} \log w_a(t) = \frac{1}{a w_a(t)^{1/p}}$, and
the constant $c_p = \frac{(1-1/p)^{p}}{p-1}$.
When $p = \infty$, each formal expression written in terms of $p$ in this proof should be interpreted as the limit of that expression as $p\to\infty$. 
For example, if $p = \infty$, $w_a(t) = \lim_{q\to\infty} (1+t/(aq))^q = e^{t/a}$ and $c_\infty = \lim_{q\to\infty} \frac{(1-1/q)^{q}}{q-1}= 0$.
For the proof of Theorem~\ref{thm:22} under the condition \eqref{eq:p1}, we introduce the energy function
\begin{talign*}%
E_k = w_a(\delta k) (f(x_{k}) - f(x^\ast)),
\end{talign*}
noting that, by the convexity of $w$ on $t \geq 0$, 
\begin{talign*}
\frac{w_a(\delta (k+1)) - w_a(\delta k)}{\delta} \leq \frac{1}{a}(1+\frac{\delta(k+1)}{ap})^{p-1} = \frac{1}{a} w_a(\delta (k+1))^{(p-1)/p}.
\end{talign*}
and hence
\begin{talign} \label{eq:discrete_w_deriv}
\frac{w_a(\delta (k+1)) - w_a(\delta k)}{\delta w_a(\delta (k+1))} \leq \frac{1}{a w_a(\delta (k+1))^{1/p}} .
\end{talign}
When \eqref{eq:p1} holds, we have
\begin{talign*}%
\frac{E_{k+1} - E_k}{\delta} 
&= \frac{w_a(\delta (k+1)) - w_a(\delta k)}{\delta}(f(x_k)-f(x^\ast)) + w_a(\delta (k+1))\frac{f(x_{k+1})-f(x_k)}{\delta} \notag \\
&\le \frac{w_a(\delta (k+1)) - w_a(\delta k)}{\delta} \langle \nabla f(x_k), x_k-x^\ast \rangle + w_a(\delta (k+1))\frac{f(x_{k+1})-f(x_k)}{\delta}  \notag\\ %
&\overset{\eqref{eq:p1}}{\le}  \frac{w_a(\delta (k+1)) - w_a(\delta k)}{\delta} \langle \nabla f(x_k), x_k-x^\ast \rangle -w_a(\delta (k+1))\|\nabla f(x_k)\|_\ast^{\frac{p}{p-1}}  \notag\\%\label{eq:E2}
&=w_a(\delta (k+1))( \frac{w_a(\delta (k+1)) - w_a(\delta k)}{\delta w_a(\delta (k+1))} \langle \nabla f(x_k), x_k-x^\ast \rangle -\|\nabla f(x_k)\|_\ast^{\frac{p}{p-1}} )  \notag\\%\label{eq:E2}
&\leq  w_a(\delta (k+1))( \frac{1}{a w_a(\delta (k+1))^{1/p}}\langle \nabla f(x_k), x_k-x^\ast \rangle -\|\nabla f(x_k)\|_\ast^{\frac{p}{p-1}} )  \notag\\
&\leq w_a(\delta (k+1))c_p\|\frac{1}{a w_a(\delta (k+1))^{1/p}}(x_k - x^\ast)\|^p \notag\\
&=c_p\|x_k - x^\ast\|^p/a^p  
\leq c_pR^p/a^p.
\end{talign*}
The first inequality uses convexity of $f$, and the second uses~\eqref{eq:p1}. 
The third inequality is an application of \eqref{eq:discrete_w_deriv}.
The fourth inequality uses the Fenchel-Young inequality  $-\|s\|^{\frac{p}{p-1}} + \langle s,u \rangle \le - \frac{p-1}{p} \|s\|^{\frac{p}{p-1}} + \langle s,u \rangle \le \frac{1}{p}\|u\|^p$
with
$s = \nabla f(x_{k})$
and
$u = \frac{1}{a w_a(\delta (k+1))^{1/p}}(x_k-x^\ast)$.
Both descent conditions~\eqref{eq:prog1} imply $\|x_k-x^\ast\| \le R$, yielding the final inequality.
Therefore, we have shown that for all $k \ge 0$, $E_{k+1}-E_k \le c_p \delta R^p/a^p.$
This implies
$E_k \le E_0 +  c_p \delta k R^p/a^p.$
Therefore
\begin{talign*}
f(x_k) - f(x^\ast) \leq \frac{f(x_0) - f(x^\ast)}{(1+\delta k/(ap))^p} + c_p \frac{R^p}{a^p}\frac{\delta k}{(1+\delta k/(ap))^p}.
\end{talign*}
Since $a > 0$ was arbitrary, we may choose $a = R\frac{(c_p\delta k)^{1/p}}{(f(x_0) - f(x^\ast))^{1/p}}$ to obtain the bound
\begin{talign*} 
f(x_k) - f(x^\ast) 
	\leq \frac{2(f(x_0) - f(x^\ast))}{\left(1+\frac{(f(x_0) - f(x^\ast))^{1/p}}{Rc_p^{1/p}p}(\delta k)^{\frac{p-1}{p}}\right)^p}
	=\textstyle O(1/(1+\frac{1}{Rp}{(\delta k)^{\frac{p-1}{p}}})^{p})
\end{talign*}
as desired. 

If, on the other hand~\eqref{eq:p2} holds,
identical reasoning yields
\begin{talign*}%
\frac{{E}_{k+1} - {E}_k}{\delta} %
&=  \frac{w_a(\delta (k+1)) - w_a(\delta k)}{\delta}(f(x_{k+1})-f(x^\ast))+ w_a(\delta k)\frac{f(x_{k+1})-f(x_k)}{\delta} \notag \\
&\le  \frac{w_a(\delta (k+1)) - w_a(\delta k)}{\delta} \langle \nabla f(x_{k+1}), x_{k+1}-x^\ast  \rangle + w_a(\delta k)\frac{f(x_{k+1})-f(x_k)}{\delta}  \notag\\ %
&\overset{\eqref{eq:p2}}{\le}  \frac{w_a(\delta (k+1)) - w_a(\delta k)}{\delta} \langle \nabla f(x_{k+1}), x_{k+1}-x^\ast \rangle -w_a(\delta k)  \|\nabla f(x_{k+1})\|_\ast^{\frac{p}{p-1}}   \notag\\%\label{eq:E2}
&=w_a(\delta k)( \frac{w_a(\delta (k+1)) - w_a(\delta k)}{\delta w_a(\delta k)} \langle \nabla f(x_{k+1}), x_{k+1}-x^\ast \rangle -\|\nabla f(x_{k+1})\|_\ast^{\frac{p}{p-1}} )  \notag\\%\label{eq:E2}
&\leq  w_a(\delta k)( \frac{w_a(\delta (k+1))}{a w_a(\delta k) w_a(\delta (k+1))^{1/p}}\langle \nabla f(x_{k+1}), x_{k+1}-x^\ast \rangle -\|\nabla f(x_{k+1})\|_\ast^{\frac{p}{p-1}} )  \notag\\
&\leq w_a(\delta k)c_p\|\frac{w_a(\delta (k+1))}{a w_a(\delta k) w_a(\delta (k+1))^{1/p}}(x_{k+1} - x^\ast)\|^p \notag\\
&=\left( \frac{w_a(\delta (k+1))}{w_a(\delta k)}\right)^{p-1}c_p\frac{R^p}{a^p}.
\end{talign*}
Now, since $w_a(\delta (k+1)) \leq w_a(\delta k) w_a(\delta)$
, we have shown
that for all $k \ge 0$, $E_{k+1}-E_k \le w_a(\delta)^{p-1} c_p \frac{R^p}{a^p}\delta.$
This implies
$E_k \le E_0 +  w_a(\delta)^{p-1} c_p \frac{R^p}{a^p}\delta k.$
Hence, we find
$$\textstyle f(x_k)-f(x^\ast) \le \frac{f(x_0) - f(x^\ast)}{(1+\delta k/(ap))^p} + w_a(\delta)^{p-1} c_p \frac{R^p}{a^p} \frac{\delta k}{(1+\delta k/(ap))^p}.$$
Since $a > 0$ was arbitrary, we may choose $a = b w_b(\delta)^{(p-1)/p}$ for $b = R\frac{(c_p\delta k)^{1/p}}{(f(x_0) - f(x^\ast))^{1/p}}$.
Since $w_b(\delta) \geq 1$, we have $b \leq a$ and hence $w_a(\delta) \leq w_b(\delta)$.
Therefore, 
\begin{talign*} 
f(x_k) - f(x^\ast) 
	\leq \frac{2(f(x_0) - f(x^\ast))}{\left(1+\frac{(f(x_0) - f(x^\ast))^{1/p}}{Rc_p^{1/p}pw_b(\delta)^{(p-1)/p}}(\delta k)^{\frac{p-1}{p}}\right)^p}
	=\textstyle O(1/(1+\frac{1}{Rp}{(\delta k)^{\frac{p-1}{p}}})^{p})
\end{talign*}
as desired.

\subsubsection{Proof of Theorem~\ref{thm:33}}
Take the energy function
$E_k = f(x_k) - f(x^\ast)$. Observe  that if~\eqref{eq:p1} holds, then we have:
\begin{talign*} 
\frac{E_{k+1} - E_k}{\delta} = \frac{f(x_{k+1}) - f(x_k)}{\delta}\overset{\eqref{eq:p1}}{\leq} -\|\nabla f(x_k)\|_{\ast}^{\frac{p}{p-1}} \overset{\eqref{eq:gdom}}{\leq} -\frac{p}{p-1}\mu^{\frac{1}{p-1}} E_k,
\end{talign*} 
or rewritten, $E_{k+1} \leq \left(1 -\frac{p}{p-1} \mu^{\frac{1}{p-1}}\delta\right)E_k$. Summing gives the bound \[E_{k+1} \leq \left(1-\frac{p}{p-1}\mu^{\frac{1}{p-1}}\delta\right)^kE_0 \leq e^{-\frac{p}{p-1}\mu^{\frac{1}{p-1}}\delta k}E_0,\] using $1+ x \leq e^x$ $\forall x \in \mathbb{R}$.
On the other hand, if~\eqref{eq:p2} holds, then a similar argument follows:
\begin{talign*} 
\frac{E_{k+1} - E_k}{\delta} = \frac{f(x_{k+1}) - f(x_k)}{\delta}\overset{\eqref{eq:p2}}{\leq} -\|\nabla f(x_{k+1})\|_{\ast}^{\frac{p}{p-1}} \overset{\eqref{eq:gdom}}{\leq} -\frac{p}{p-1}\mu^{\frac{1}{p-1}} E_{k+1},
\end{talign*} 
or rewritten, $E_{k+1} \leq \left(1 + \frac{p}{p-1}\mu^{\frac{1}{p-1}}\delta\right)^{-1}E_k$. Summing gives the bound 
\begin{talign*} E_{k+1} \leq \left(1+ \frac{p}{p-1}\mu^{\frac{1}{p-1}}\delta\right)^{-k}E_0 \leq e^{-\frac{p}{p-1}\mu^{\frac{1}{p-1}}\delta k}E_0.\end{talign*}

\subsection{Examples of descent methods}
\label{sec:examples}
We now provide detailed demonstration that the examples provided are descent algorithms.
\subsubsection{Higher-order gradient descent}

Let $ \tilde p = p-1+\nu$.
The optimality condition for the HGD algorithm~\eqref{eq:hgd} is
\begin{talign}\label{Eq:UpdateProof1}
\sum_{i=1}^{p-1} \frac{1}{(i-1)!} \nabla^i f(x_k) \, (x_{k+1}-x_k)^{i-1} + \frac{1}{\eta} \|x_{k+1}-x_k\|^{\tilde p-2} \,B (x_{k+1}-x_k) = 0.
\end{talign}
Since $\nabla^{p-1} f$ is $L$-Lipschitz, we have the following error bound on the $(p-2)$-nd order Taylor expansion of $\nabla f$:
\begin{talign}\label{Eq:UpdateProof2}
\left\|\nabla f(x_{k+1}) - \sum_{i=1}^{p-1} \frac{1}{(i-1)!} \nabla^i f(x_k) \, (x_{k+1}-x_k)^{i-1}\right\|_* \le \frac{L}{(p-2)!} \|x_{k+1}-x_k\|^{p-2+\nu}.
\end{talign}
Substituting~\eqref{Eq:UpdateProof1} to~\eqref{Eq:UpdateProof2} and writing $r_k = \|x_{k+1}-x_k\|$, we obtain
\begin{talign}\label{Eq:UpdateProof2a}
\left\|\nabla f(x_{k+1}) + \frac{r_k^{\tilde p-2}}{\eta} \,B (x_{k+1}-x_k)\right\|_* \,\le\, \frac{L}{(p-2)!} r_k^{\tilde p-1}.
\end{talign}
Squaring both sides, expanding, and rearranging the terms, we get the inequality
\begin{talign}\label{Eq:UpdateProof3}
\langle \nabla f(x_{k+1}), x_k-x_{k+1} \rangle
\,\ge\, \frac{\eta}{2r_k^{\tilde p-2}} \|\nabla f(x_{k+1})\|_*^2 + \frac{\eta r_k^{\tilde p}}{2}\left(\frac{1}{\eta^2} - \frac{L^2}{(p-2)!^2}\right).
\end{talign}
If $p=2$, then the first term in~\eqref{Eq:UpdateProof3} already implies the desired bound below. Now assume $p \ge 3$. The right-hand side of~\eqref{Eq:UpdateProof3} is of the form $A/r^{\tilde p-2} + Br^{\tilde p}$, which is a convex function of $r > 0$ and minimized by $r^* = \left\{\frac{(\tilde p-2)}{\tilde p} \frac{A}{B} \right\}^{\frac{1}{2\tilde p-2}}$, yielding a minimum value of
\begin{talign*} 
\frac{A}{(r^*)^{\tilde p-2}} + B(r^*)^{\tilde p}
\,=\, A^{\frac{p}{2 \tilde p-2}} B^{\frac{\tilde p-2}{2\tilde p-2}} \left[\left(\frac{\tilde p}{\tilde p-2}\right)^{\frac{\tilde p-2}{2\tilde p-2}} + \left(\frac{\tilde p-2}{\tilde p}\right)^{\frac{\tilde p}{\tilde p-2}}\right]
\,\ge\, A^{\frac{p}{2\tilde p-2}} B^{\frac{\tilde p-2}{2\tilde p-2}}.
\end{talign*}
Substituting the values $A = \frac{\eta}{2} \|\nabla f(x_{k+1})\|_*^2$ and $B = \frac{\eta}{2} (\frac{1}{\eta^2} - \frac{L^2}{(p-2)!^2})$ from~\eqref{Eq:UpdateProof3}, we obtain
\begin{talign*} 
\langle \nabla f(x_{k+1}), x_k-x_{k+1} \rangle
\,&\ge\, \frac{\eta}{2} \left(\frac{1}{\eta^2} - \frac{L^2}{(p-2)!^2}\right)^{\frac{\tilde p-2}{2\tilde p-2}} \|\nabla f(x_{k+1})\|_*^{\frac{\tilde p}{\tilde p-1}}.
\end{talign*}
Finally, using the inequality $f(x_k) - f(x_{k+1}) \ge \langle \nabla f(x_{k+1}), x_k-x_{k+1} \rangle$ by the convexity of $f$ yields the progress bound
\begin{talign*} 
f(x_{k+1}) - f(x_k) &\leq  -\frac{\eta^{\frac{1}{\tilde p-1}}}{2}\left(1 - \frac{(L\eta)^2}{(p-2)!^2}\right)^{\frac{\tilde p-2}{2\tilde p-2}} \|\nabla f(x_{k+1})\|_*^{\frac{\tilde p}{\tilde p-1}}\\&\leq -\frac{\eta^{\frac{1}{\tilde p-1}}}{2^{\frac{2\tilde p-3}{\tilde p-1}}}  \|\nabla f(x_{k+1})\|_*^{\frac{\tilde p}{\tilde p-1}}
\end{talign*}
where the least inequality uses the fact that $\eta \leq \frac{\sqrt{3}(p-2)!}{2L}$.
\subsubsection{Proximal method}

The optimality condition for the proximal method is
$$\nabla^2 h(x_k)^{-1}\nabla f(x_{k+1}) + \frac{\|x_{k+1}-x_k\|_{x_k}^{p-2}}{\eta}(x_{k+1}-x_k) = 0,$$
which implies $\|x_{k+1}-x_k\|_{x_k} = \eta^{\frac{1}{p-1}} \|\nabla f(x_{k+1})\|_{x_k,\ast}^{\frac{1}{p-1}}$, using the shorthand $\|v\|_{x_k,\ast}  = \sqrt{ \langle v, \nabla h(x_k)^{-1} v\rangle }$.
From the definition of $x_{k+1}$, we have $f(x_{k+1}) + \frac{1}{p\eta}\|x_{k+1}-x_k\|_{x_k}^p \le f(x_k)$.
Rearranging gives
$$\textstyle f(x_k) - f(x_{k+1}) \ge \frac{1}{p\eta}\|x_{k+1}-x_k\|_{x_k}^p = \frac{\eta^{\frac{1}{p-1}}}{p} \|\nabla f(x_{k+1})\|_{\ast,x_k}^{\frac{p}{p-1}}\ge  \frac{m^{\frac{p}{p-1}}\eta^{\frac{1}{p-1}}}{p} \|\nabla f(x_{k+1})\|_\ast^{\frac{p}{p-1}}$$
as desired.

\subsubsection{Natural gradient descent}

Since $\nabla^2 f \preceq LB$, we have the bound
$$\textstyle f(x_{k+1}) \le f(x_k) + \langle \nabla f(x_k), x_{k+1}-x_k \rangle + \frac{L}{2}\|x_{k+1}-x_k\|^2.$$
Plugging in the NGD update~\eqref{eq:ngd} gives
$$\textstyle f(x_{k+1}) \le f(x_k) - \eta \langle \nabla f(x_k), (\nabla^2 h(x_k))^{-1} \nabla f(x_k) \rangle + \frac{L \eta^2}{2} \langle \nabla f(x_k), B(\nabla^2 h(x_k))^{-2} \nabla f(x_k) \rangle.$$
Since $m B \preceq \nabla^2 h \preceq MB$, we have $\frac{1}{M}B^{-1} \preceq (\nabla^2 h)^{-1} \preceq \frac{1}{m} B^{-1}$, so
\begin{talign*} 
f(x_{k+1}) &\le f(x_k) - \frac{\eta}{M} \|\nabla f(x_k)\|_\ast^2 + \frac{L \eta^2}{2 m^2} \|\nabla f(x_k)\|_\ast^2 \\
&= f(x_k) - \eta \left(\frac{1}{M} - \frac{L \eta}{2m^2}\right) \|\nabla f(x_k)\|_\ast^2 \\
&\le f(x_k)  - \frac{\eta}{2M} \|\nabla f(x_k)\|_\ast^2
\end{talign*}
where in the last step we have used the inequality $\eta \le \frac{m^2}{ML}$.

\subsubsection{Mirror descent}
Plugging the variational condition $\nabla h(x_{k+1}) - \nabla h(x_k) = -\eta \nabla f(x_k)$ into the smoothness bound on $f$, as well as using the property $mB  \preceq \nabla^2 h$  we have%
\begin{talign*} 
f(x_{k+1}) - f(x_k) &\leq \langle \nabla f(x_k) , x_{k+1} - x_k\rangle + \frac{L}{2}\|x_{k+1} - x_k\|^2\\%&\le- \frac{1}{\eta}\langle \nabla h(x_{k+1}) - \nabla h(x_k), x_{k+1}-x_k \rangle + \frac{L}{2}\|x_{k+1}-x_k\|^2\\
&\leq- \frac{1}{\eta}\langle \nabla h(x_{k+1}) - \nabla h(x_k), x_{k+1}-x_k \rangle + \frac{L}{2m^2}\|\nabla h(x_{k+1})-\nabla h(x_k)\|_\ast^2
\end{talign*}
Given $h$ is $M$-smooth, 
$- \frac{1}{\eta}\langle \nabla h(x_{k+1}) - \nabla h(x_k), x_{k+1}-x_k \rangle \leq -\frac{1}{\eta M}\|\nabla h(x_{k+1}) - \nabla h(x_k)\|_\ast^2$
(\cite[(2.1.8)]{Nesterov04}) and therefore,
\begin{talign*} 
f(x_{k+1}) - f(x_k) \le -\left(\frac{1}{\eta M}-  \frac{L}{2m^2}  \right)\|\nabla h(x_{k+1}) - \nabla h(x_k)\|_\ast^2 &\leq - \eta \left(\frac{1}{M} - \frac{L \eta}{2m^2}\right) \|\nabla f(x_k)\|_\ast^2\\
& \le- \frac{\eta}{2M} \|\nabla f(x_k)\|_\ast^2
\end{talign*}
where in the last step we have used the inequality $\eta \le \frac{m^2}{ML}$. %
\subsubsection{Proximal Bregman Method}
The optimality condition for the proximal method is
$\eta \nabla f(x_{k+1}) =   \nabla h(x_{k+1}) - \nabla h(x_k)$, which implies $\eta^2\|\nabla f(x_{k+1})\|_\ast^2 = \| \nabla h(x_{k+1}) - \nabla h(x_k)\|_\ast^2 \leq M^2\|x_{k+1} - x_k\|^2$. From the definition of $x_{k+1}$, we have $ f(x_{k+1}) + \frac{1}{\eta} D_h(x_{k+1}, x_k) \le f(x_k) $. Rearranging gives
$$\textstyle  f(x_{k+1}) - f(x_k) \le - \frac{1}{\eta}D_h(x_{k+1},x_k) \leq   -\frac{m}{2\eta}\|x_{k+1} - x_k\|^2 \leq -\frac{m\eta}{2M^2}\|\nabla f(x_{k+1})\|_\ast^2$$
as desired.

\subsection{Rescaled Gradient Descent}
\label{App:RGD}
\paragraph{Proof of Lemma~\ref{lem:1}} 
We show rescaled gradient descent satisfies progress bound~\eqref{eq:prog1} with $\delta = \eta^{\frac{1}{p-1}}/2$ when $f$ is strongly smooth. Since $\|\nabla ^p f(x)\|\leq L_p$, we have the Taylor expansion bound,
\begin{talign*}
f(x_{k+1}) - f(x_k) &\leq \langle \nabla f(x_k), x_{k+1} - x_k\rangle +  \sum_{m=2}^{p-1}\frac{1}{m!} \nabla^m f(x_k)(x_{k+1} - x_k)^m + \frac{L_p}{p!}\|x_{k+1} - x_k\|^p\\
&\overset{\eqref{eq:rgd}}{=} - \eta^{\frac{1}{p-1}}\left(1 - \frac{\eta L_p}{p!}\right)\|\nabla f(x_k)\|_\ast^{\frac{p}{p-1}} +  \sum_{m=2}^{p-1}\frac{\eta^{\frac{m}{p-1}}}{m!} \frac{\nabla^m f(x_k)(\nabla f(x_k))^m}{\|\nabla f(x_k)\|_\ast^{\frac{m(p-2)}{p-1}}}\\
&\overset{\eqref{eq:ss}}{\leq}   - \eta^{\frac{1}{p-1}}\left(1 - \frac{\eta L_p}{p!}\right)\|\nabla f(x_k)\|_\ast^{\frac{p}{p-1}} +  \sum_{m=2}^{p-1}\frac{\eta^{\frac{m}{p-1}}}{m!} L_m\|\nabla f(x_k)\|_\ast^{m + \frac{p-m}{p-1} - \frac{m(p-2)}{p-1}}\\
& =   - \eta^{\frac{1}{p-1}}\left(1 - \frac{\eta L_p}{p!}\right)\|\nabla f(x_k)\|_\ast^{\frac{p}{p-1}} +  \sum_{m=2}^{p-1}\frac{\eta^{\frac{m}{p-1}}}{m!}L_m \|\nabla f(x_k)\|_\ast^{\frac{p}{p-1}}\\
 &= - \eta^{\frac{1}{p-1}}\left(1 -\sum_{m=2}^{p}\frac{\eta^{\frac{m-1}{p-1}} L_m}{m!}\right)\|\nabla f(x_k)\|_\ast^{\frac{p}{p-1}}. 
\end{talign*}
The second line follows from the rescaled gradient update~\eqref{eq:rgd} and the third follows from our strongly smoothness Assumption (def ~\ref{ass:ss2}). Since $\eta<1$ we can further bound
\begin{talign*}
f(x_{k+1}) - f(x_k) &\leq- \eta^{\frac{1}{p-1}}\left(1 -\eta^{\frac{1}{p-1}}\sum_{m=2}^{p}\frac{ L_m}{m!}\right)\|\nabla f(x_k)\|_\ast^{\frac{p}{p-1}}. 
\end{talign*}
Our step-size condition~\eqref{eq:step} implies $ 1-\eta^{\frac{1}{p-1}}\sum_{m=2}^{p}\frac{ L_m}{m!} \geq \frac{1}{2}$, which  yields the desired bound~\eqref{eq:prog1} with $\delta = \eta^{\frac{1}{p-1}}/2$.

\subsection{Gradient Descent vs.\ Rescaled Gradient Descent}
\label{App:GDSlow}

\paragraph{Proof of Lemma~\ref{lem:1}} 
We have $f'(x) = \mathrm{sign}(x) |x|^{p-1}$, so $|f'(x)|^{\frac{p-2}{p-1}} = |x|^{p-2}$.

The rescaled gradient descent of order $p$ with step size $\epsilon = \eta^{\frac{1}{p-1}}$ is
\begin{align*}
x_{k+1} = x_k - \epsilon \frac{f'(x_k)}{|f'(x_k)|^{\frac{p-2}{p-1}}} = x_k - \epsilon\frac{\mathrm{sign}(x_k) |x_k|^{p-1}}{|x_k|^{p-2}} = (1-\epsilon) x_k.
\end{align*}
Therefore, if $0 < \epsilon < 1$, then $x_k = (1-\epsilon)^k x_0$, and thus $f(x_k) = (1-\epsilon)^{pk} f(x_0)$ converges to $0$ at an exponential rate $\Theta((1-\epsilon)^{pk})$.

The gradient descent with step size $\epsilon = \eta^{\frac{1}{p-1}}$ for $f$ is
\begin{align*}
x_{k+1} = x_k - \epsilon f'(x_k) = x_k - \epsilon \, \mathrm{sign}(x_k) |x_k|^{p-1} = (1-\epsilon |x_k|^{p-2}) x_k.
\end{align*}
Note that if $0 < \epsilon < |x_k|^{-(p-2)}$, then $x_{k+1}$ has the same sign as $x_k$ with smaller magnitude.
In particular, if $0 < x_0 < \epsilon^{-\frac{1}{p-2}}$, then $x_k > x_{k+1} > 0$ for all $k > 0$,
and gradient descent simplifies to $x_{k+1} = (1-\epsilon x_k^{p-2})x_k$.
Assume we start with $0 < x_0 \le (2\epsilon)^{-\frac{1}{p-2}}$, so $\frac{x_k}{x_{k+1}} = (1-\epsilon x_k^{p-2})^{-1} \le  (1-\epsilon x_0^{p-2})^{-1} \le 2$.
Then by Jensen's inequality applied to the convex function $x \mapsto x^{-(p-2)}$, we have
$x_{k+1}^{-(p-2)} - x_k^{-(p-2)} \le \frac{-(p-2)}{x_{k+1}^{p-1}} (x_k-x_{k+1}) =(p-2) \epsilon \frac{x_k^{p-1}}{x_{k+1}^{p-1}}  \le (p-2) 2^{p-1} \epsilon$.
This implies $x_k \ge (x_0^{-(p-2)} + (p-2) 2^{p-1} \epsilon k)^{-\frac{1}{p-2}} = \Omega((\epsilon k)^{-\frac{1}{p-2}})$, and thus $f(x_k) \ge \Omega((\epsilon k)^{-\frac{p}{p-2}})$ converges to $0$ at a polynomial rate.

\medskip

\subsubsection{Gradient Flow vs.\ Rescaled Gradient Flow}

We also discuss how the behavior in discrete time above matches the behavior in continuous time.
The rescaled gradient flow of order $p$ for $f$ is
\begin{align*}
\dot X_t = -\frac{f'(X_t)}{|f'(X_t)|^{\frac{p-2}{p-1}}} = -\frac{\mathrm{sign}(X_t) |X_t|^{p-1}}{|X_t|^{p-2}} = -X_t
\end{align*}
so $X_t = e^{-t} X_0$, and thus $f(X_t) = e^{-pt} f(X_0)$ converges to $0$ at an exponential rate $\Theta(e^{-pt})$.

The gradient flow (which is rescaled gradient flow of order $2$) for $f$ is
\begin{align*}
\dot X_t = -f'(X_t) = -\mathrm{sign}(X_t) |X_t|^{p-1}
\end{align*}
Without loss of generality assume $X_0 > 0$, so $X_t > 0$ for all $t > 0$.
Then gradient flow simplifies to $\dot X_t = -X_t^{p-1}$, or $\frac{d}{dt} X_t^{-(p-2)} = -(p-2) \dot X_t \, X_t^{-(p-1)} = p-2$, so $X_t = (X_0^{-(p-2)} + (p-2)t )^{-\frac{1}{p-2}}$, and thus $f(X_t) = \Theta(t^{-\frac{p}{p-2}})$ converges to $0$ at a polynomial rate.

More generally, the rescaled gradient flow of order $q$ ($q > 1, q \neq p$) for $f$ is
\begin{align*}
\dot X_t = -\frac{f'(X_t)}{|f'(X_t)|^{\frac{q-2}{q-1}}} = -\frac{\mathrm{sign}(X_t) |X_t|^{p-1}}{|X_t|^{\frac{(q-2)(p-1)}{q-1}}} 
= - \mathrm{sign}(X_t) |X_t|^{\frac{p-1}{q-1}}
\end{align*}
Assume $X_0 > 0$, so $X_t > 0$ for all $t > 0$.
Rescaled gradient flow simplifies to $\dot X_t = -X_t^{\frac{p-1}{q-1}}$, or $\frac{d}{dt} X_t^{-\frac{p-q}{q-1}} = \frac{p-q}{q-1} $, so $X_t = (X_0^{-\frac{p-q}{q-1}} + (\frac{p-q}{q-1})t )^{-\frac{q-1}{p-q}}$, and $f(X_t) = \Theta(t^{-\frac{p(q-1)}{p-q}})$.
Note that if $1 < q < p$, then $f(X_t)$ converges to $0$ at a polynomial rate, which becomes faster as $q \to p$.
At $q=p$, the convergence rate becomes exponential, as we see for rescaled gradient flow above.
However, for $q > p$, $f(X_t)$ diverges to $\infty$.
Thus, the best order to use is $q=p$, but it is better to underestimate $p$.

\section{Accelerating Descent Algorithms}
\label{app:accel}
The energy function
\begin{talign}\label{eq:lyap}
E_k = D_h(x^\ast, z_k) + A_k(f(y_k) - f(x^\ast)) ,
\end{talign}
will be used to analyze all the accelerated methods introduced in this paper. 

\subsection{Proof of Proposition~\ref{prop:8}}
\label{app:acc}
Take energy (Lyapunov) function~\eqref{eq:lyap} Set $A_k = C\delta^p k^{(p)}$ where $ k^{(p)} = k(k+1)\cdots(k+p-1)$ is the rising factorial. Denote $\alpha_k := \frac{A_{k+1} - A_k}{\delta} = Cp\delta^{p-1}(k+1)^{(p-1)}$ and $\tau_k := \frac{\alpha_k}{A_{k+1}} = \frac{k}{\delta(k+p)}$.
\paragraph{Algorithm~\eqref{eq:nest1}:} Using~\eqref{eq:lyap} we compute
\begin{talign}\label{eq:pro1}
 \frac{E_{k+1} - E_k}{\delta} = \frac{D_h(x^\ast, z_{k+1}) - D_h(x^\ast, z_k)}{\delta}+ \frac{A_{k+1}}{\delta}(f(y_{k+1}) - f(x^\ast)) - \frac{A_{k}}{\delta} (f(y_{k}) - f(x^\ast)).
\end{talign}
We bound the first part,
 \begin{talign}
\frac{D_h(x^\ast, z_{k+1}) - D_h(x^\ast, z_k)}{\delta} &= - \left\langle \frac{\nabla h(z_{k+1}) - \nabla h(z_k)}{\delta}, x^\ast - z_{k+1}\right\rangle - \frac{1}{\delta} D_h(z_{k+1}, z_k)\notag\\
&\overset{\eqref{eq:up2}}{=}\alpha_k\langle \nabla f(x_{k}), x^\ast - z_k\rangle +\alpha_k \langle \nabla f(x_{k}), z_k - z_{k+1}\rangle - \frac{1}{\delta} D_h(z_{k+1}, z_k)\notag\\
& \leq \alpha_k\langle \nabla f(x_{k}), x^\ast - z_k\rangle   + (\delta/m)^{\frac{1}{p-1}}\alpha_k^{\frac{p}{p-1}}\|\nabla f(x_{k})\|_\ast^{\frac{p}{p-1}}, \label{eq:pro2} %
\end{talign}
where the inequality follows from the $m$-uniform convexity of $h$ of order $p$ and the Fenchel-Young inequality $\langle s,h\rangle + \frac{1}{p}\|h\|^p \geq -\frac{p}{p-1} \|s\|_\ast^{\frac{p}{p-1}}\leq -\|s\|_\ast^{\frac{p}{p-1}}$, with $h =(m/\delta)^{\frac{1}{p}}(z_{k+1} - z_k)$ and $s = (\delta/m)^{\frac{1}{p}}\alpha_k^{\frac{p}{p-1}} \nabla f(x_k)$. Plugging in update~\eqref{eq:up1}, %
\begin{talign}
 \alpha_k\langle \nabla f(x_{k}), x^\ast - z_k\rangle &= \alpha_k\langle \nabla f(x_{k}), x^\ast - y_k\rangle + \frac{A_{k+1}}{\delta}\langle \nabla f(x_{k}), y_k - x_{k}\rangle \notag\\
 &= \alpha_k\langle \nabla f(x_{k}), x^\ast - x_k\rangle + \frac{A_{k}}{\delta} \langle \nabla f(x_k), y_k - x_k\rangle\notag\\
& \leq - \left(\frac{A_{k+1}}{\delta}(f(y_{k+1}) - f(x^\ast)) -\frac{A_{k}}{\delta} (f(y_{k}) - f(x^\ast))\right) \notag\\&\quad+  A_{k+1}\frac{f(y_{k+1}) - f(x_k)}{\delta}\notag\\% \label{eq:pro3}\\
&\leq-  \left(\frac{A_{k+1}}{\delta}(f(y_{k+1}) - f(x^\ast)) -\frac{A_{k}}{\delta} (f(y_{k}) - f(x^\ast))\right) \notag\\&\quad-  A_{k+1}\delta^{\frac{1}{p-1}}\|\nabla f(x_k)\|_\ast^{\frac{p}{p-1}}.\label{eq:pro3}
\end{talign}
The first inequality follows from the convexity of $f$ and rearranging terms. The second inequality uses the progress condition assumed for the sequence $y_{k+1}$. Combining~\eqref{eq:pro1} with \eqref{eq:pro2} and \eqref{eq:pro3} we have, 
\begin{talign*} 
\frac{E_{k+1} - E_k}{\delta} \leq \left((\delta/m)^{\frac{1}{p-1}}(Cp\delta^{p-1}(k+1)^{(p-1)})^{\frac{p}{p-1}}  -C\delta^{\frac{1}{p-1}}\delta^p(k+1)^{(p)}\right)\|\nabla f(x_{k})\|_\ast^{\frac{p}{p-1}}. 
\end{talign*}
Given $((k+1)^{(p-1)})^{\frac{p}{p-1}}/(k+1)^{(p)}\leq 1$, it suffices that
$C \leq 1/mp^{p} $ to ensure $\frac{E_{k+1} - E_k}{\delta} \leq 0$. Summing the Lyapunov function gives the convergence rate $f(y_k) - f(x^\ast) = O(1/A_k)= O(1/(\delta k)^p)$.
\paragraph{Algorithm~\eqref{eq:nest2}:}
Using~\eqref{eq:lyap} with the same parameter choices as algorithm~\eqref{eq:nest1}, we have
  \begin{talign}
\frac{D_h(x^\ast, z_{k+1}) - D_h(x^\ast, z_k)}{\delta} & \leq \alpha_k\langle \nabla f(y_{k+1}), x^\ast - z_k\rangle  +(\delta/m)^{\frac{1}{p-1}}\alpha_k^{\frac{p}{p-1}}\|\nabla f(y_{k+1})\|_\ast^{\frac{p}{p-1}},  \label{eq:pro22}
\end{talign}
where the first part uses the same steps as~\eqref{eq:pro2} except update~\eqref{eq:up22} is used instead of~\eqref{eq:up2}.
Plugging in update~\eqref{eq:up11} yields the following, 
\begin{talign}
\alpha_k\langle \nabla f(y_{k+1}), x^\ast - z_k\rangle &= \alpha_k\langle \nabla f(y_{k+1}), x^\ast - y_{k+1}\rangle +\frac{A_{k+1}}{\delta}\langle \nabla f(y_{k+1}), y_{k+1} - z_{k}\rangle \notag\\
 &\overset{\eqref{eq:up11}}{=} \alpha_k\langle \nabla f(y_{k+1}), x^\ast - y_{k+1}\rangle + \frac{A_k}{\delta} \langle \nabla f(y_{k+1}),y_{k} - y_{k+1} \rangle\notag\\
  & \quad+ \frac{A_{k+1}}{\delta}\langle \nabla f(y_{k+1}), y_{k+1} - x_k\rangle\notag\\
& \leq - \left(\frac{A_{k+1}}{\delta}(f(y_{k+1}) - f(x^\ast)) -\frac{A_k}{\delta} (f(y_{k}) - f(x^\ast))\right)\notag \\
&\quad+\frac{A_{k+1}}{\delta}\langle \nabla f(y_{k+1}), y_{k+1} - x_k\rangle \notag\\ %
& \leq - \left(\frac{A_{k+1}}{\delta}(f(y_{k+1}) - f(x^\ast)) -\frac{A_{k}}{\delta}(f(y_{k}) - f(x^\ast))\right)\notag\\ &\quad- A_{k+1}\delta^{\frac{1}{p-1}}\|\nabla f(y_{k+1})\|_\ast^{\frac{p}{p-1}}. \label{eq:pro33}
\end{talign}
The first inequality follows from the convexity of $f$ and rearranging terms. The second inequality uses the progress condition assumed for the sequence $y_{k+1}$. Combining~\eqref{eq:pro1} with~\eqref{eq:pro22} \eqref{eq:pro33}, we have
\begin{talign*} 
\frac{E_{k+1} - E_k}{\delta} %
&\leq-\delta^{\frac{1}{p-1}} C(k+1)^{(p)}\|\nabla f(y_{k+1})\|_\ast^{\frac{p}{p-1}} +(\delta /m)^{\frac{1}{p-1}}(Cp(k+1)^{(p-1)})^{\frac{p}{p-1}}\|\nabla f(y_{k+1})\|_\ast^{\frac{p}{p-1}}.
\end{talign*}
For $\frac{E_{k+1} - E_k}{\delta} \leq 0$  it suffices that
$C \leq1/mp^{p} $. Summing the Lyapunov function gives the convergence rate $f(y_k) - f(x^\ast) = O(1/A_k)= O(1/(\delta k)^p)$.

\subsection{Restarting Scheme}
\label{app:restart}

When $f$ is {\em strongly smooth} and $\mu$-gradient dominated, we define the restarting scheme (similar to~\cite[(B.1.2)]{acceleration}), which proceeds by running~\ref{eq:argd} 
for some number of iterations at each step,
  \begin{align}\label{eq:rest}
  \hat x_k = (\text{the output } y_c \text{ of running Algorithm~\ref{eq:argd} for $c$ iterations with input $x_0 = \hat x_{k-c}$).}
  \end{align}
\begin{theorem}
Assume $f$ is convex and strongly smooth of order $1<p<\infty$  with constants $0<L_1,\dots,L_p<\infty$ {\em and} $f$ is $\mu$-gradient dominated of order $p$. Suppose $\eta$ satisfies~\eqref{eq:step}. Let $\hat x_k$ be the output of running the restarting scheme~\eqref{eq:rest} for $k/c$ times with $c = 2p/\kappa^\frac{1}{p}$ where $\kappa = \mu\delta^p = \mu \eta$. Finally, let $y_k$ be the output of running the rescaled gradient descent update one step from $\hat x_k$. The composite scheme  satisfies the convergence rate upper bound: 
$f( y_k) - f(x^\ast) =O(  \exp({-\frac{1}{2p}\mu^{\frac{1}{p}}\delta k})$
\end{theorem}

Take $h(x) = \frac{2^{p-2}}{p}\|x - x_0\|^p$ which is $1$-uniformly convex of order $p$. Running $k$ iterations of either algorithm~\eqref{eq:nest1} or~\eqref{eq:nest2} results in the convergence bound,
\begin{talign}
\frac{\mu}{p}\|\hat x_k - x^\ast \|^p \leq f(\hat x_k) - f(x^\ast) \leq \frac{2^{p-2}p^{p-1}\|\hat{x}_{k-c} - x^\ast\|^p}{\delta^{p}k^p}&\leq \frac{  2^{p-2} p^{p-1}\|\hat{x}_{k-c} - x^\ast\|^p}{(\delta c)^p} \notag\\%\notag\\
&\leq \frac{\mu}{p e}\|\hat{x}_{k-c} - x^\ast\|^p.\label{eq:refer}
\end{talign}
where the last inequality follows from the choice $c = 2p/\kappa^{\frac{1}{p}}$. Thus an execution of~\eqref{eq:rest} for $c$ iterations of the accelerated method reduces the distance to optimum by a factor of at least $1/e$. Iterating~\eqref{eq:refer}, we obtain $\frac{1}{p}\|\hat x_k - x^\ast \|^p \leq e^{-k/c} \frac{1}{p}\|\hat x_0 - x^\ast \|^p$. Using the descent property for both methods, $E_{k+1} \leq \delta 2p^{p-1}\|x_k - x^\ast\|^p $~\eqref{eq:p1} and $E_{k+1} \leq \delta 2p^{p-1}\|x_{k+1} - x^\ast\|^p $~\eqref{eq:p2}, implies that 
\begin{talign*} 
f(\hat y_{k}) - f(x^\ast) \leq \delta 2p^{p-1}e^{\frac{-\kappa^{\frac{1}{p}}k}{2p}}\|x_0 - x^\ast\|^p = O\left(e^{\frac{-\kappa^{\frac{1}{p}}k}{2p}}\right).
\end{talign*}

\subsection{Proof of Proposition~\ref{prop:MS}}
\label{App:MSAcc}
We analyze the following sequence of iterates
\begin{subequations}\label{eq:nest3}
\begin{talign}
x_{k} &= \delta\tau_k z_k + (1 - \delta\tau_k) y_k\label{eq:up111}\\
z_{k+1} &= \arg \min_{z}\left\{ \alpha_k \langle \nabla f(y_{k+1}), z\rangle +  \frac{1}{\delta}D_h(z, z_k)\right\}\label{eq:up222},
\end{talign}
where the update for $(\lambda_{k+1},y_{k+1})$ satisfies the descent conditions
\begin{talign}
a \leq \frac{\lambda_{k+1}}{\delta^{\frac{3p-2}{2}}}\|y_{k+1} - x_k\|^{p-2}&\leq b,\label{eq:descent1}\\
\|y_{k+1} - x_{k} + \frac{\lambda_{k+1}}{m} \nabla f(y_{k+1})\|& \leq \frac{1}{2} \|y_{k+1} - x_{k}\|,\label{eq:descent2}
\end{talign}
\end{subequations}
and the following identifications $\alpha_k = \frac{ A_{k+1} - A_k}{\delta}$,  $\tau_k = \frac{\alpha_k}{A_{k+1}}$,  and $\lambda_{k+1} = \frac{\alpha_k^2}{\delta^2 A_{k+1}}$ hold. Assume $h$ is $m$-strongly convex.

Taking energy function~\eqref{eq:lyap}, we compute 
\begin{talign*} 
\frac{E_{k+1} - E_k}{\delta} &= \frac{A_{k+1}}{\delta}( f(y_{k+1}) - f(x^\ast)) - \frac{A_k}{\delta} (f(y_k) - f(x^\ast))\\
&\quad   -\left\langle \frac{\nabla h(z_{k+1}) - \nabla h(z_k)}{\delta}, x^\ast - z_{k+1}\right\rangle- D_h(z_{k+1} ,z_{k})\\
&\overset{\eqref{eq:up222}}{\leq} \alpha_k (f(y_{k+1}) - f(x^\ast)) + \frac{A_k}{\delta} (f(y_{k+1}) - f(y_k)) +\alpha_k\langle \nabla f(y_{k+1}), x^\ast - z_{k+1}\rangle \\
&\quad- \frac{m}{2\delta}\|z_{k} - z_{k+1}\|^2\\
&\leq \alpha_k \langle \nabla f(y_{k+1}), y_{k+1} - z_{k+1}\rangle +\frac{A_k}{\delta}\langle \nabla f(y_{k+1}), y_{k} - y_{k+1}\rangle -  \frac{m}{2\delta}\|z_{k} - z_{k+1}\|^2.
\end{talign*}
where the first inequality follows from the strong convexity of $h$ and the last inequality follows from the convexity of $f$. Denote $x = \delta\tau_k  z_{k+1} + (1 - \delta \tau_k) y_k$. 
Starting from the preceding line, we have,
\begin{talign*} 
\frac{E_{k+1} - E_k}{\delta} &\leq\frac{A_{k+1}}{\delta} \left\langle \nabla f(y_{k+1}), y_{k+1} - x \right\rangle  -\frac{m}{2\delta}\|z_{k} - \frac{1}{\delta\tau_k} x + \frac{1-\delta\tau_k}{\delta\tau_k}y_{k}\|^2\\
&= \frac{A_{k+1}}{\delta} \left\langle \nabla f(y_{k+1}), y_{k+1}- x\right\rangle  -\frac{1}{2(\delta \tau_k)^2}\frac{m}{\delta}\|\delta\tau_kz_{k} +(1-\delta\tau_k)y_k  -x \|^2\\
&\overset{\eqref{eq:up111}}{=} \frac{A_{k+1}}{\delta} \left\langle \nabla f(y_{k+1}),  y_{k+1} - x\right\rangle  -\frac{1}{2(\delta \tau_k)^2}\frac{m}{\delta}\|x_{k} -x \|^2\\
&\leq \max_{x \in \mathcal{X}}\left\{ \frac{A_{k+1}}{\delta} \left\langle \nabla f(y_{k+1}), y_{k+1} - x\right\rangle  -\frac{1}{2(\delta \tau_k)^2}\frac{m}{\delta}\|x_{k} -x \|^2\right\}.
\end{talign*}
Plugging in the solution, which satisfies $x = x_{k} - \frac{\delta^2}{m}\frac{\alpha_k^2}{A_{k+1}} \nabla f(y_{k+1})$, and noting $\lambda_{k+1} = \frac{\delta^2\alpha_k^2}{A_{k+1}}$ we obtain
\begin{talign}
\frac{E_{k+1} - E_k}{\delta} &\leq  \frac{A_{k+1}}{\lambda_{k+1}}\frac{m}{2\delta}\left( \|y_{k+1} - x_k +\frac{\lambda_{k+1} }{m}\nabla f(y_{k+1})\|^2 - \| y_{k+1}  - x_k\|^2\right) \notag\\
&\overset{\eqref{eq:descent2}}{\leq } -  \frac{A_{k+1}}{\lambda_{k+1}\delta}\frac{m}{4}\|y_{k+1} - x_k\|^2. \label{eq:conv_bound}
\end{talign}
This is the same bound as~\cite[(3.12)]{Svaiter} with $\sigma = 0$.

Rearranging the last inequality and summing over $k$, we have 
\begin{talign}\label{eq:last}
 \sum_{i=0}^k \frac{A_{i}}{\lambda_i} \frac{m}{4}\|y_{i+1} - x_i\|^2 \leq E_{k+1} +  \sum_{i=0}^k \frac{A_{i}}{\lambda_i} \frac{m}{4}\|y_{i+1} - x_i\|^2\leq E_0 = D_h(x^\ast,x_0),
\end{talign}
where the last equality comes from taking $A_0 = 0$. 

Notice that summing over our bound~\eqref{eq:conv_bound} gives us the rate
\begin{align*}
f(y_k) - f(x^\ast) \leq \frac{E_0}{A_k}.
\end{align*}
Now we use the second bound~\eqref{eq:descent1} to establish $A_k = O(k^{\frac{3p-2}{2}})$. This follows from arguments identical to the those given by~\cite[p.6-7]{Gasnikov} and~\cite[p.6-8]{Bubeck}. Denote $a_1 = a \delta^{\frac{3p-2}{2}}$. Observe that
\begin{talign} \label{eq:last3}
\sum_{i=0}^k \frac{A_{i}}{\lambda_{i}^{\frac{p}{p-2}}} a_1^{\frac{2}{p-2}} \overset{\eqref{eq:descent1}}{\leq}  \sum_{i=0}^k \frac{A_{i}}{\lambda_{i}^{1 + \frac{2}{p-2}}}\left(
\lambda_i\|y_{i+1} - x_i\|^{p-2}\right)^{\frac{2}{p-2}} \leq\sum_{i=0}^k \frac{A_{i}}{\lambda_{i}}\|y_{i+1} - x_i\|^2 \overset{\eqref{eq:last}}{\leq} 4 E_0/m.
\end{talign}
Denote $c_1 = a_1^{-\frac{2}{p-2}} 4E_0/m = (a\delta^{\frac{3p-2}{2}})^{-\frac{2}{p-2}} E_04/m$. Using the previous line, we have
\begin{talign}\label{eq:bubeck_in} 
A_k \geq \frac{1}{4}\left(\sum_{i=1}^k \sqrt{\lambda_i} \right)^2\geq \frac{1}{4}c_1^{-\frac{p-2}{p}}\left(\sum_{i=1}^{k}A_i^{\frac{p-2}{3p-2}}\right)^{\frac{3p-2}{p}},
\end{talign}
where the first inequality follows from definition of $\alpha_k$ (see~\cite[Lem 2.6]{Bubeck}) and the second inequality uses reverse Holders (see~\cite[p.7-8]{Bubeck}).
Specifically, we have
$$ \alpha_k = \frac{\lambda_k + \sqrt{\lambda_k^2 + 4\lambda_k A_{k-1}}}{2} \geq \frac{\lambda_k }{2} + \sqrt{\lambda_k A_{k-1}} \geq \left( \frac{\lambda_k }{2} + \sqrt{ A_{k-1}}
\right)^2 - A_{k-1},$$
and $\alpha_k^2 = \lambda_k A_k$ which allows us to conclude the first inequality. For the second inequality, we use reverse Holder (i.e. $\|fg\|_1 \geq \|f\|_{\frac{1}{q}}\|g\|_{-\frac{1}{q-1}}$ for $q \geq 1$) with $q = 1 + \frac{p-2}{2p} = \frac{3p -2}{2p}$ so that $-\frac{1}{q-1} = \frac{2p}{p-2}$, we have 
\begin{talign}\label{eq:last1} \sum_{i=0}^k \sqrt{\lambda_i} = \sum_{i=0}^k A_i^{\frac{p-2}{2p}}\left(\frac{A_i}{\lambda_i^{\frac{p}{p-2}}}\right)^{-\frac{p-2}{2p}} \geq \left( \sum_{i=0}^k A_i^{\frac{p-2}{3p-2}}\right)^{\frac{3p-2}{2p}}\left( \sum_{i=0}^k \frac{A_i}{\lambda_i^{\frac{p}{p-2}}}\right)^{-\frac{p-2}{2p}}.
\end{talign}
Equation \eqref{eq:bubeck_in} follows from combining~\eqref{eq:last1} with~\eqref{eq:last3}.

To end our proof, we use the elementary fact~\cite[Lem  3.4]{Bubeck} that for a positive sequence $B_j$ such that $B_k^\alpha \geq c_2 \sum_{i=1}^k B_j$, we have 
\begin{talign*} 
B_k \geq \left(\frac{\alpha -1 }{\alpha} c_2 k\right)^{\frac{1}{\alpha-1}}
\end{talign*}
with the identificatons $\alpha = \frac{p}{p-2}$, $B_k  = A_k^{\frac{p-2}{3p-2}}$ and $c_2 = \frac{c_1^{-\frac{p-2}{3p-2}}}{4^{\frac{p}{3p-2}}  } $. 
Subsequently,
\[
A_k \geq \left(\frac{2c_2 k}{p}\right)^{\frac{3p-2}{2}} = \Theta\left((\delta k)^{\frac{3p-2}{2}}E_0^{-\frac{p-2}{2}}\right),\] as desired. Picking up the constants, we have the bound
\begin{talign*} 
f(y_k) - f(x^\ast) \leq \frac{E_0}{A_k} = \frac{c_3D_h(x^\ast, x_0)^{\frac{p}{2}}}{(\delta k)^{\frac{3p-2}{2}}},
\end{talign*}
where $c_3^{-1} = a(2/p)^{\frac{3p-2}{2}} (4/m)^{-\frac{p-2}{2}}$. %

\subsection{Restarting Scheme}
\label{app:restart2}
When $f$ is {\em strongly smooth} and $\mu$-gradient dominated, we define the restarting scheme (similar to~\eqref{eq:rest}), which proceeds by running Algorithm~\ref{alg:argd2}  
for some number of iterations at each step,
  \begin{align}\label{eq:rest2}
  \hat x_k = (\text{the output } y_c \text{ of running Algorithm~\ref{alg:argd2} for $c$ iterations with input $x_0 = \hat x_{k-c}$).}
  \end{align}

We summarize the behavior of the restarting scheme in the following theorem:

\begin{theorem}
Assume $f$ is convex and $s$-strongly smooth of order $1 < p <\infty$ with constants $0<L_1, \dots, L_p<\infty$ and $f$ is $\mu$-gradient dominated of order $p$. Take $h(x) = \frac{1}{2}\|x\|^2$. Let $\hat x_k$ be the output of running the restarting scheme~\eqref{eq:rest2} for $k/c$ times with $c =(p^3/2)^{\frac{p}{3p-2}}(e/3\kappa)^{\frac{2}{3p-2}}$ where $\kappa = \mu \delta^{\frac{3p-2}{2}}= \mu \eta$. Finally, let $y_k$ be the output of running the rescaled gradient descent update one step from $\hat x_k$. Then we have the convergence rate upper bound: 
\begin{align*} 
f( y_k) - f(x^\ast) = O\left( \exp\left(-c_1\mu^{\frac{2}{3p-2}}\delta k\right)\right),
\end{align*}
where $c_1 =  (3/e)^{\frac{2}{3p-2}} (2/p^3)^{\frac{p}{3p-2}}$. 
\end{theorem}

\label{app:restart2}
Take $h(x) = \frac{1}{2}\|x\|^2$ which is $1$-strongly convex. Running $k$ iterations of algorithm~\eqref{eq:nest3} results in the convergence bound
\begin{talign}\label{eq:ref}
\frac{\mu}{p}\|\hat x_k - x^\ast \|^p \leq f(\hat x_k) - f(x^\ast) \leq \frac{\frac{c_3}{2}\|\hat{x}_{k-c} - x^\ast\|^p}{(\delta k)^{\frac{3p-2}{2}}}\leq \frac{ \frac{c_3}{2}\|\hat{x}_{k-c} - x^\ast\|^p}{(\delta c)^{\frac{3p-2}{2}}}\leq \frac{\mu}{p e}\|\hat{x}_{k-c} - x^\ast\|^p,
\end{talign}
where the last inequality follows from the choice $c =(c_3pe/2\kappa)^{\frac{2}{3p-2}}$ where $\kappa = \delta^{\frac{3p-2}{2}}\mu$. Thus an execution of~\eqref{eq:rest2} for $c$ iterations of the accelerated method reduces the distance to optimum by a factor of at least $1/e$. Iterating~\eqref{eq:ref}, we obtain $\frac{1}{p}\|\hat x_k - x^\ast \|^p \leq e^{-k/c} \frac{1}{p}\|\hat x_0 - x^\ast \|^p$. Here, we require that the update from $x_k$ to $y_{k+1}$ be a descent algorithm. Using the descent property for both methods $E_{k+1} \leq \delta 2p^{p-1}\|x_k - x^\ast\|^p $~\eqref{eq:p1} and $E_{k+1} \leq \delta 2p^{p-1}\|x_{k+1} - x^\ast\|^p $~\eqref{eq:p2} implies that 
\begin{talign*} 
f(\hat y_{k}) - f(x^\ast) \leq \delta 2p^{p-1}e^{-c_4\mu^{\frac{2}{3p-2}}\delta k}\|x_0 - x^\ast\|^p = O\left(e^{-c_4\mu^{\frac{2}{3p-2}}\delta k}\right),
\end{talign*}
where $c_4 = (c_3pe/2)^{-\frac{2}{3p-2}}$.

\subsection{Proof of Theorem~\ref{thm:M-ARGD}}
\label{app:thmM-ARGD}
We show under the strong smoothness, rescaled gradient descent with line search condition~\eqref{eq:descent1} satisfies~\eqref{eq:descent2}. We summarize in the following Lemma. 
\begin{lemma}
Under the above assumptions, if $\eta^{\frac{1}{p-1}} \le \min\{\frac{2}{5p},1/(2 \sum_{m=2}^p \frac{L_m}{m!})\}$ and $\lambda_{k+1}$ is such that
\begin{talign}\label{Eq:Cond1}
\frac{3}{4} \le \frac{\lambda_{k+1} \|x_{k+1}-x_k\|^{p-2}}{\eta} \le \frac{5}{4},
\end{talign}
then rescaled gradient descent~\eqref{eq:rgd} satisfies
\begin{talign}\label{Eq:Want1}
\|x_{k+1} - x_k + \lambda_{k+1} \nabla f(x_{k+1})\| \le \frac{1}{2} \|x_{k+1}-x_k\|.
\end{talign}
\end{lemma}
 Note, we can write~\eqref{Eq:Cond1} as
\begin{talign}\label{Eq:Cond}
\frac{3}{4} \frac{\eta^{\frac{1}{p-1}}}{\|\nabla f(x_k)\|^{\frac{p-2}{p-1}}} \le \lambda_{k+1} \le \frac{5}{4} \frac{\eta^{\frac{1}{p-1}}}{\|\nabla f(x_k)\|^{\frac{p-2}{p-1}}}.
\end{talign}
Plugging in the RGD update~\eqref{eq:rgd} to~\eqref{Eq:Want1}, what we wish to show is that
\begin{talign}\label{Eq:Want2}
\left\|\lambda_{k+1} \nabla f(x_{k+1}) - \frac{\eta^{\frac{1}{p-1}}}{\|\nabla f(x_k)\|^{\frac{p-2}{p-1}}} \nabla f(x_k) \right\| \le \frac{\eta^{\frac{1}{p-1}}}{2} \|\nabla f(x_k)\|^{\frac{1}{p-1}}.
\end{talign}

Since $\|\nabla^p f(x)\| \le L_p$, we have the following Taylor expansion of $\nabla f$:
$$\textstyle \nabla f(x_{k+1}) = \nabla f(x_k) + \sum_{m=2}^{p-1} \frac{1}{(m-1)!} (\nabla^m f(x_k))(x_{k+1}-x_k)^{m-1} + R_k$$
where $R_k$ is the remainder term which can be bounded as
$$\textstyle \|R_k\| \le \frac{L_p}{(p-1)!} \|x_{k+1}-x_k\|^{p-1} = \frac{L_p}{(p-1)!} \eta \|\nabla f(x_k)\|.$$
Furthermore, by strong smoothness assumption, for $m = 2,\dots,p-1$ we have
\begin{talign*} 
\|(\nabla^m f(x_k))(x_{k+1}-x_k)^{m-1}\| 
&= \eta^{\frac{m}{p-1}} \frac{|(\nabla^m f(x_k))(\nabla f(x_k))^{m-1}|}{\|\nabla f(x_k)\|^{\frac{(m-1)(p-2)}{p-1}}} \\
&\le \eta^{\frac{m}{p-1}} \frac{L_m \|\nabla f(x_k)\|^{m-1+\frac{p-m}{p-1}}}{\|\nabla f(x_k)\|^{\frac{(m-1)(p-2)}{p-1}}} \\
&= \eta^{\frac{m}{p-1}} L_m \|\nabla f(x_k)\|.
\end{talign*}

By plugging in the bounds above to the left-hand side of~\eqref{Eq:Want2}, we get
\begin{talign*} 
&\left\|\lambda_{k+1} \nabla f(x_{k+1}) - \frac{\eta^{\frac{1}{p-1}}}{\|\nabla f(x_k)\|^{\frac{p-2}{p-1}}} \nabla f(x_k) \right\| \\
&= \left\| \left(\lambda_{k+1} - \frac{\eta^{\frac{1}{p-1}}}{\|\nabla f(x_k)\|^{\frac{p-2}{p-1}}}\right) \nabla f(x_k) + \lambda_{k+1} \sum_{m=2}^{p-1} \frac{1}{(m-1)!} (\nabla^m f(x_k))(x_{k+1}-x_k)^{m-1} + \lambda_{k+1} R_k \right\| \\
&\le \left| \lambda_{k+1} - \frac{\eta^{\frac{1}{p-1}}}{\|\nabla f(x_k)\|^{\frac{p-2}{p-1}}}\right| \|\nabla f(x_k)\| + \lambda_{k+1} \sum_{m=2}^{p-1} \frac{1}{(m-1)!} \|(\nabla^m f(x_k))(x_{k+1}-x_k)^{m-1}\| + \lambda_{k+1} \|R_k\| \\
&\le \left| \lambda_{k+1} - \frac{\eta^{\frac{1}{p-1}}}{\|\nabla f(x_k)\|^{\frac{p-2}{p-1}}}\right| \|\nabla f(x_k)\| + \lambda_{k+1} \sum_{m=2}^{p-1} \frac{1}{(m-1)!} \eta^{\frac{m}{p-1}} L_m \|\nabla f(x_k)\|_\ast + \lambda_{k+1} \frac{L_p}{(p-1)!} \eta \|\nabla f(x_k)\|\\
&= \left(\left| \lambda_{k+1} - \frac{\eta^{\frac{1}{p-1}}}{\|\nabla f(x_k)\|^{\frac{p-2}{p-1}}}\right| + \lambda_{k+1} \sum_{m=2}^{p-1} \frac{\eta^{\frac{m}{p-1}} L_m}{(m-1)!} + \lambda_{k+1} \frac{L_p}{(p-1)!} \eta \right) \|\nabla f(x_k)\| \\
&= \left(\left| \lambda_{k+1} - \frac{\eta^{\frac{1}{p-1}}}{\|\nabla f(x_k)\|^{\frac{p-2}{p-1}}}\right| + \lambda_{k+1} \sum_{m=2}^{p} \frac{\eta^{\frac{m}{p-1}} m L_m}{m!} \right) \|\nabla f(x_k)\| \\
&\le \left(\left| \lambda_{k+1} - \frac{\eta^{\frac{1}{p-1}}}{\|\nabla f(x_k)\|^{\frac{p-2}{p-1}}}\right| + \lambda_{k+1} \eta^{\frac{2}{p-1}} p \sum_{m=2}^{p} \frac{L_m}{m!} \right) \|\nabla f(x_k)\| \\
&\le \left(\left| \lambda_{k+1} - \frac{\eta^{\frac{1}{p-1}}}{\|\nabla f(x_k)\|^{\frac{p-2}{p-1}}}\right| + \lambda_{k+1} \frac{\eta^{\frac{1}{p-1}} p}{2} \right) \|\nabla f(x_k)\|
\end{talign*}
where in the last step we have used that $\eta^{\frac{1}{p-1}} \le 1/(2 \sum_{m=2}^p \frac{L_m}{m!})$.

Therefore, from the above, we see that if
\begin{talign}\label{Eq:Want2a}
\left| \lambda_{k+1} - \frac{\eta^{\frac{1}{p-1}}}{\|\nabla f(x_k)\|^{\frac{p-2}{p-1}}}\right| \le \frac{\eta^{\frac{1}{p-1}}}{4\|\nabla f(x_k)\|^{\frac{p-2}{p-1}}}
\end{talign}
and
\begin{talign}\label{Eq:Want2b}
\lambda_{k+1} \frac{\eta^{\frac{1}{p-1}} p}{2} \le \frac{\eta^{\frac{1}{p-1}}}{4\|\nabla f(x_k)\|^{\frac{p-2}{p-1}}},
\end{talign}
then the desired relation~\eqref{Eq:Want2} holds.
The first condition~\eqref{Eq:Want2a} is equivalent to
$$\textstyle \frac{3}{4} \frac{\eta^{\frac{1}{p-1}}}{\|\nabla f(x_k)\|^{\frac{p-2}{p-1}}} \le \lambda_{k+1} \le \frac{5}{4} \frac{\eta^{\frac{1}{p-1}}}{\|\nabla f(x_k)\|^{\frac{p-2}{p-1}}}$$
which is precisely the requirement~\eqref{Eq:Cond}, whereas the second condition~\eqref{Eq:Want2b} is equivalent to
$$\textstyle \lambda_{k+1} \le \frac{1}{2p \|\nabla f(x_k)\|^{\frac{p-2}{p-1}}}.$$
Note that if $\eta^{\frac{1}{p-1}} \le \frac{2}{5p}$, then the last condition above is automatically satisfied if the right-hand side of the former condition~\eqref{Eq:Cond} holds.
Therefore, we have shown that the condition~\eqref{Eq:Cond} implies the desired relation~\eqref{Eq:Want2}, or equivalently~\eqref{Eq:Want1}. A simple continuity argument, similar to~\cite[Lem 3.2]{Bubeck} ensures the existence of pair $(\lambda_k, y_k)$ that satisfies~\eqref{Eq:Cond1} and~\eqref{Eq:Want1} simultaneously.

\subsection{Proximal method}
\label{thm:P-ARGD}
Given $x_k \in \R^n$ and $\eta > 0$, let $x_{k+1}$ be the proximal update~\eqref{eq:prox}, which satisfies 
\begin{talign}\label{Eq:Proximal}
x_{k+1} = x_k - \eta^{\frac{1}{p-1}} \frac{\nabla f(x_{k+1})}{\|\nabla f(x_{k+1})\|^{\frac{p-2}{p-1}}}.
\end{talign}

\begin{lemma}
If $\lambda_{k+1}$ is such that
 \begin{talign}\label{eq:this}
\frac{1}{2} \le \frac{\lambda_{k+1} \|x_{k+1}-x_k\|^{p-2}}{\epsilon} \le \frac{3}{2},
\end{talign}
then
\begin{talign}\label{Eq:Want3}
\|x_{k+1} - x_k + \lambda_{k+1} \nabla f(x_{k+1})\| \le \frac{1}{2} \|x_{k+1}-x_k\|.
\end{talign}
\end{lemma}
Note~\eqref{eq:this} is equivalent to the condition
 \begin{talign}\label{Eq:Cond2}
\frac{1}{2} \frac{\eta^{\frac{1}{p-1}}}{\|\nabla f(x_{k+1})\|^{\frac{p-2}{p-1}}} \le \lambda_{k+1} \le \frac{3}{2} \frac{\eta^{\frac{1}{p-1}}}{\|\nabla f(x_{k+1})\|^{\frac{p-2}{p-1}}}.
\end{talign}
Plugging in the proximal update~\eqref{Eq:Proximal} to~\eqref{Eq:Want3}, what we wish to show is that
$$\textstyle \left\|\lambda_{k+1} \nabla f(x_{k+1}) - \frac{\eta^{\frac{1}{p-1}}}{\|\nabla f(x_{k+1})\|^{\frac{p-2}{p-1}}} \nabla f(x_{k+1}) \right\| \le \frac{\eta^{\frac{1}{p-1}}}{2} \|\nabla f(x_{k+1})\|^{\frac{1}{p-1}}.$$
Equivalently, we wish to show that
$$\textstyle \left| \lambda_{k+1} - \frac{\eta^{\frac{1}{p-1}}}{\|\nabla f(x_{k+1})\|^{\frac{p-2}{p-1}}}\right| \le \frac{\epsilon^{\frac{1}{p-1}}}{2\|\nabla f(x_{k+1})\|^{\frac{p-2}{p-1}}},$$
which is exactly condition~\eqref{Eq:Cond2}. 
Subsequently, we can write the Monteiro-Svaiter-style accelerated proximal method as the following sequence of updates,
\begin{algorithm}[H]
\caption{Monteiro-Svaiter-style accelerated proximal method}
\label{alg:argd5}
\begin{algorithmic}[1]
\Require{$f$ is differentiable and $h$ is $1$-strongly convex}
	\State Set $x_0 = z_0 =0$, $A_0 = 0$, $\delta^{\frac{3p-2}{2}} = \eta$, $\eta >0$\\
	{\bf for} $k = 1, \dots, K$ {\bf do}
	\begin{subequations}%
		\State Choose $\lambda_{k+1}$ (e.g. by line search) such that 
	$\frac{1}{2} \leq \frac{\lambda_{k+1} \|y_{k+1} -  x_k\|^{p-2}}{\eta} \leq \frac{3}{2}$,
	where 
	\begin{talign*} 
	y_{k+1} = \arg \min_{x \in \mathcal{X}}\left\{f(x) + \frac{1}{\eta p}\|x - x_k\|^p\right\},
	\end{talign*}
	\end{subequations}
	and $\alpha_{k} =\frac{ \lambda_{k+1} + \sqrt{\lambda_{k+1} + 4A_k\lambda_{k+1}}}{2\delta}$, $A_{k+1} = \delta \alpha_{k}  + A_k$, $\tau_k = \frac{\alpha_k}{A_{k+1}}$ (so that $\lambda_{k+1} = \frac{\delta^2 \alpha_k^2}{A_{k+1}}$) and %
	$$\textstyle x_k = \delta \tau_k  z_k + (1-\delta \tau_k) y_k.$$

	\State Update $z_{k+1} = \arg \min_{z\in \mathcal{X}}\left\{ \alpha_k \langle \nabla f(y_{k+1}), z\rangle +   \frac{1}{\delta}D_h(z, z_k)\right\}$
\State \Return $y_K$.
\end{algorithmic}
\end{algorithm}

\section{Examples and Numerical Experiments}
\label{supp:exp}
\subsection{Comparison to Runge-Kutta}
In~\cite{Runge} the following gradient lower bound assumption is made
\begin{definition}
$f$ satisfies the {\em gradient lower bound} of order $p \ge 2$ if for all $m = 1,\dots,p-1$,
$$\textstyle f(x)-f(x^\ast) \ge \frac{1}{C_m} \|\nabla^m f(x)\|^{\frac{p}{p-m}} ~~~ \forall \, x \in \R^n$$
for some constants $0 < C_1, \dots, C_{p-1} < \infty$.
\end{definition}
Notice that when $p=2$, this is equivalent to $s$-strong smoothness, which is the general smoothness condition on the gradient. However, for $p>2$ we can show that it is slightly weaker than strong smoothness. We summarize in the following Lemma:
\begin{lemma}
If $f$ is strongly smooth of order $p$ with constants $L_m$, then $f$ satisfies the gradient lower bound of order $p$ with constants $C_m = 4(\sum_{i=2}^p \frac{L_i}{i!})L_m^{\frac{p}{p-m}}$.
\end{lemma}
Let $\eta = 1/(2 \sum_{m=2}^p \frac{L_m}{m!})^{p-1}$ as in~\eqref{ass:ss2}.
Then with $x_k = x$ and $x_{k+1} = x - \eta^{\frac{1}{p-1}} \nabla f(x)/\|\nabla f(x)\|^{\frac{p-2}{p-1}}$, by Lemma~\ref{lem:1} we have
$$\textstyle f(x^\ast) \le f(x_{k+1}) \le f(x) -\frac{\eta^{\frac{1}{p-1}}}{2} \|\nabla f(x)\|^{\frac{p}{p-1}}
= f(x) - \frac{1}{4 \sum_{m=2}^p \frac{L_m}{m!}} \|\nabla f(x)\|^{\frac{p}{p-1}}.$$
Rearranging gives the desired claim:
\begin{talign*}
f(x) - f(x^\ast) \ge \frac{1}{4 \sum_{m=2}^p \frac{L_m}{m!}} \|\nabla f(x)\|^{\frac{p}{p-1}}.
\end{talign*}

\subsection{Examples}
\label{app:example}
We provide details on the examples presented in the main text. 
\label{app:example}
\subsection{$\ell_p$ loss}
Let
$$f(x) = \frac{1}{p} \|x\|_p^p = \frac{1}{p} \sum_{i=1}^d |x_i|^p = \frac{1}{p} \sum_{i=1}^d \text{sgn}(x_i)^p x_i^p.$$
The gradient $\nabla f(x)$ has entries
$$(\nabla f(x))_i = \text{sgn}(x_i)^p x_i^{p-1}.$$
The norm of the gradient is
$$\|\nabla f(x)\| = \left(\sum_{i=1}^d x_i^{2p-2}\right)^{\frac{1}{2}} = \|x\|_{2p-2}^{p-1}.$$
Therefore, for $m \ge 2$,
$$\|\nabla f(x)\|^{\frac{p-m}{p-1}} = \|x\|_{2p-2}^{p-m} = \left(\sum_{i=1}^d x_i^{2p-2}\right)^{\frac{p-m}{2p-2}}.$$

For $m \ge 2$, the $m$-th derivative $\nabla^m f(x)$ has nonzero entries only on the diagonal:
$$(\nabla^m f(x))_{i,\dots,i} = (p-1)\cdots(p-m+1) \text{sgn}(x_i)^p x_i^{p-m}.$$
Then for any unit vector $v \in \R^d$,
$$(\nabla^m f(x))(v^m) = (p-1)\cdots(p-m+1) \sum_{i=1}^d \text{sgn}(x_i)^p x_i^{p-m} v_i^m.$$
By H\"older's inequality with $q = \frac{2p-2}{p-m}$ and $r = \frac{2p-2}{p+m-2}$, so $\frac{1}{q} + \frac{1}{r} = 1$, we have
\begin{align*}
|(\nabla^m f(x))(v^m)| &= (p-1)\cdots(p-m+1) \left|\sum_{i=1}^d \text{sgn}(x_i)^p x_i^{p-m} v_i^m \right| \\
&\le (p-1)\cdots(p-m+1) \left(\sum_{i=1}^d |\text{sgn}(x_i) x_i^{p-m}|^{\frac{2p-2}{p-m}}\right)^{\frac{p-m}{2p-2}} \left(\sum_{i=1}^d |v_i^m|^{\frac{2p-2}{p+m-2}}\right)^{\frac{p+m-2}{2p-2}} \\
&= (p-1)\cdots(p-m+1) \, \|x\|_{2p-2}^{p-m} \left(\sum_{i=1}^d |v_i|^{\frac{2m(p-1)}{p+m-2}}\right)^{\frac{p+m-2}{2p-2}}.
\end{align*}
Note that $\frac{m(p-1)}{p+m-2} = 1 + \frac{(m-1)(p-2)}{p+m-2} \ge 1$.
Then using $\sum_{i=1}^d c_i^q \le (\sum_{i=1}^d c_i)^q$ for $c_i \ge 0$, $q \ge 1$, we can write
$$\sum_{i=1}^d |v_i|^{\frac{2m(p-1)}{p+m-2}} \le \left(\sum_{i=1}^d v_i^2 \right)^{\frac{m(p-1)}{p+m-2}} = \|v\|_2^{\frac{2m(p-1)}{p+m-2}} = 1$$
since we assumed $v$ is a unit norm vector, so $\|v\|_2 = 1$.
Plugging this to the bound above, we obtain
\begin{align*}
|(\nabla^m f(x))(v^m)| &\le (p-1)\cdots(p-m+1) \, \|x\|_{2p-2}^{p-m}  \\
&= (p-1)\cdots(p-m+1) \, \|\nabla f(x)\|^{\frac{p-m}{p-1}}.
\end{align*}
Taking the supremum over unit vectors $v \in \R^d$, we conclude that
$$\|\nabla^m f(x)\| \le (p-1)\cdots(p-m+1) \|\nabla f(x)\|^{\frac{p-m}{p-1}}.$$
This shows that $f$ is strongly smooth of order $p$ with constants
$$L_m = (p-1)\cdots(p-m+1).$$
\subsection{Logistic loss}
We show the logistic loss of strongly smooth of order $p = \infty$.
We have
$$\nabla f(x) = -\frac{w}{1+e^{-w^\top x}}$$
and
$$\|\nabla f(x)\| = \frac{\|w\|}{1+e^{-w^\top x}}.$$
By induction we can see that
$$\nabla^m f(x) = -\frac{(m-1)! w^{\otimes m}}{(1+e^{-w^\top x})^m}$$
so that
$$\|\nabla^m f(x)\| = \sup_{\|v\| = 1} |(\nabla^m f(x))(v^m)| = \frac{(m-1)! \|w\|^{m}}{(1+e^{-w^\top x})^m}.$$
Then
$$\frac{\|\nabla^m f(x)\|}{\|\nabla f(x)\|} = \frac{(m-1)! \|w\|^{m-1}}{(1+e^{-w^\top x})^{m-1}} \le (m-1)! \|w\|^{m-1}.$$
This shows that $f(x) = \log(1+e^{-w^\top x})$ satisfies the strong smoothness condition with $p = \infty$ with constant
$$L_m = (m-1)! \|w\|^{m-1}.$$
\subsection{GLM loss}
Consider the generalized linear model loss function $f(x) = \frac{1}{2} (y - \phi(x^\top w))^2$ for $\phi(r) = 1/(1+e^{-r}) \in (0,1)$, $y \in \{0,1\}$, and $w \in \R^d$.
Introduce the shorthand $b = 1-2y \in \{1, -1\}$, and note that 
\begin{talign*}
\phi(r) - y &= b\phi(b r), \\
\phi'(r) &= e^{-r}/(1+e^{-r})^2 = \phi(r) \phi(-r) = \phi'(-r) \in (0,1/4],\\
\phi'(r) /\phi(r) &= \phi(-r), \\
\phi''(r) 
&= \phi'(r) \phi(-r) - \phi(r) \phi'(-r) 
= \phi'(r) (\phi(-r) - \phi(r)) \in [-1/(6\sqrt{3}),1/(6\sqrt{3})], \\
\phi''(r) /\phi'(r) &= \phi(-r) - \phi(r), \qtext{and}\\
\phi'''(r) &= \phi''(r) (\phi(-r) - \phi(r)) - 2\phi'(r)^2 \in [-1/2, 0]
\end{talign*}
To simplify the presentation, we will fix $x$ and let $z = x^\top w$.
With this notation in place we have
\begin{talign*}
f(x) &= \frac{1}{2} \phi(b z)^2,\\
\nabla f(x) 
&= b \phi(b z) \phi'(b z) w
,\\
\nabla^2 f(x) 
&= (\phi'(b z)^2 +  \phi(b z) \phi''(b z)) ww^\top, \qtext{and} \\
\nabla^3 f(x) 
&= b(3\phi'(b z) \phi''(b z)+\phi(b z) \phi'''(b z)) w^{\otimes 3}.
\end{talign*}
Since $\phi(r)\phi'(r) \in (0,1)$, 
we have, for any $a \in [0,1]$
\begin{talign*}
\frac{\|\nabla^2 f(x)\|}{\|\nabla f(x)\|^a }  
	&= \frac{|\phi'(b z)^2 +  \phi(b z) \phi''(b z)|}{|\phi(b z) \phi'(b z)|^a} \|w\|^{2-a}
	\leq \frac{|\phi'(b z)^2 +  \phi(b z) \phi''(b z)|}{|\phi(b z) \phi'(b z)|} \|w\|^{2-a}\\
	&= |2\phi(-b z) - \phi(bz)| \|w\|^{2-a}
	\leq 2 \|w\|^{2-a}.
\end{talign*}
Moreover,
\begin{talign*}
\\|\nabla^3 f(x)\|
	&= |3\phi'(b z) \phi''(b z)+\phi(b z) \phi'''(b z)|  \|w\|^{3}
	\leq (\sqrt{3}/24 + 1/2) \|w\|^{3}.
\end{talign*}
Therefore, $f$ is s-strongly smooth of order $p = 3$ with $L_2 = 2 \|w\|^{1.5}$ and $L_3 = (\sqrt{3}/24 + 1/2) \|w\|^{3}$.

\section{Additional Results}

\subsection{Coordinate Descent Methods}
\label{app:Cor}
At each iteration, a randomized coordinate method samples a coordinate direction $i \in \{1, \dots, d\}$ uniformly at random and performs an update along that coordinate direction. Denote $\nabla_{i_k} f = e_{i_k} e_{i_k}^\top \nabla f(x)$ where $e_{i}$ is the $i$-th basis vector.
\begin{definition}
An algorithm $x_{k+1} = \mathcal{A}(x_k)$ is a {\bf coordinate descent algorithm of order $1<p\leq\infty$}, if for some constant $0<\delta<\infty$, it almost surely satisfies
\begin{align}\label{eq:corprog1}
\frac{f(x_{k+1}) - f(x_k)}{\delta} \leq   -\|\nabla_{i_k} f(x_k)\|_\ast^{\frac{p}{p-1}}.
\end{align}
\end{definition}
For coordinate descent methods of order $p$, it is possible to obtain non-asymptotic guarantees for non-convex, convex and gradient dominated functions. We summarize in the following theorems. 
\begin{theorem} \label{thm:111}
Suppose an algorithm satisfies~\eqref{eq:corprog1} for some $0<\delta<\infty$ and $1<p\leq\infty$ and $f$ is differentiable. Then the algorithm also satisfies:
\begin{align}\label{eq:b111}
\min_{0\leq s\leq k}\mathbb{E} \|\nabla_{i_s} f(x_s)\|_\ast \leq (E_0/(\delta k))^{\frac{p-1}{p}} = O(1/\delta k).
\end{align}
\end{theorem}
\begin{theorem}\label{thm:222}
Suppose an algorithm satisfies~\eqref{eq:corprog1} for some $0<\delta<\infty$ and $1<p\leq \infty$ and $f$ is differentiable and convex with $R = \sup_{x:f(x)\leq f(x_0)} \|x - x^\ast \| < \infty$.  
 Then the algorithm satisfies: 
\begin{align}\label{eq:b222}
\mathbb{E}[f(x_k)] - f(x^\ast) = 
\begin{cases}
O\left(1/\left(1+\frac{1}{Rp}{(\delta k)^{\frac{p-1}{p}}}\right)^{p}\right)&\text{ if }\, p < \infty\\
O\left(e^{-\delta k/R}\right)&\text{ if }\, p = \infty
 \end{cases}.
\end{align}
\end{theorem}

\begin{theorem}\label{thm:333}
Suppose an algorithm satisfies~\eqref{eq:prog1} for some $0<\delta<\infty$ and $1<p\leq\infty$, and $f$ is differentiable and $\mu$-gradient dominated of order $p$. Then the algorithm satisfies:
\begin{align}\label{eq:b333}
\mathbb{E}[f(x_k)] - f(x^\ast) = O\left( e^{-\frac{1}{d}\frac{p}{p-1}\mu^{\frac {1}{p-1}}\delta k}\right).
\end{align}
\end{theorem}

\subsubsection{Proof of Theorem~\ref{thm:111}}
\begin{talign*} 
 \delta k \mathbb{E}\min_{0\leq s \leq k} \|\nabla_s f(x_{s})\|_{\ast}^{\frac{p}{p-1}} \leq \mathbb{E}\sum_{s={0}}^k \|\nabla_s f(x_{s})\|_{\ast}^{\frac{p}{p-1}} \delta\leq f(x_0) - \mathbb{E}f(x_k) \leq f(x_0)
\end{talign*}
Rearranging the inequality yields the result in Theorem~\ref{thm:111}.\\

\subsubsection{Proof of Theorem~\ref{thm:222}}
For the proof of Theorem~\ref{thm:222} under the condition \eqref{eq:corprog1}, we use the energy function
\begin{talign*}%
E_k = w_a(\delta k) (f(x_{k}) - f(x^\ast)),
\end{talign*}
When \eqref{eq:corprog1} holds, we have
\begin{talign*}%
\frac{E_{k+1} - E_k}{\delta} 
&= \frac{w_a(\delta (k+1)) - w_a(\delta k)}{\delta}(f(x_k)-f(x^\ast)) + w_a(\delta (k+1))\frac{f(x_{k+1})-f(x_k)}{\delta} \notag \\
&\le \frac{w_a(\delta (k+1)) - w_a(\delta k)}{\delta} \langle \nabla f(x_k), x_k-x^\ast \rangle + w_a(\delta (k+1))\frac{f(x_{k+1})-f(x_k)}{\delta}  \notag\\ %
&\overset{\eqref{eq:corprog1}}{\le}  \frac{w_a(\delta (k+1)) - w_a(\delta k)}{\delta} \langle \nabla f(x_k), x_k-x^\ast \rangle -w_a(\delta (k+1))\|\nabla_{i_k} f(x_k)\|_\ast^{\frac{p}{p-1}}  \notag\\%\label{eq:E2}
&=w_a(\delta (k+1))( \frac{w_a(\delta (k+1)) - w_a(\delta k)}{\delta w_a(\delta (k+1))} \langle \nabla f(x_k), x_k-x^\ast \rangle -\|\nabla_{i_k} f(x_k)\|_\ast^{\frac{p}{p-1}} )  \notag\\%\label{eq:E2}
&\leq  w_a(\delta (k+1))( \frac{1}{a w_a(\delta (k+1))^{1/p}}\langle \nabla f(x_k), x_k-x^\ast \rangle -\|\nabla_{i_k} f(x_k)\|_\ast^{\frac{p}{p-1}} )  \notag\\
&=  w_a(\delta (k+1))( \frac{1}{a w_a(\delta (k+1))^{1/p}}\langle \nabla_{i_k} f(x_k), x_k-x^\ast \rangle -\|\nabla_{i_k} f(x_k)\|_\ast^{\frac{p}{p-1}} )  +\xi_k\notag\\
&\leq w_a(\delta (k+1))c_p\|\frac{1}{a w_a(\delta (k+1))^{1/p}}(x_k - x^\ast)\|^p + \xi_k \notag\\
&=c_p\|x_k - x^\ast\|^p/a^p  + \xi_k
\leq c_pR^p/a^p + \xi_k.
\end{talign*}
Here, the martingale $\xi_k :=\frac{ w_a(\delta (k+1))}{a w_a(\delta (k+1))^{1/p}}\langle\nabla f(x_k)- \nabla_{i_k} f(x_k), x_k-x^\ast \rangle $.
The first inequality uses convexity of $f$, and the second uses~\eqref{eq:p1}. 
The third inequality is an application of \eqref{eq:discrete_w_deriv}.
The fourth inequality uses the Fenchel-Young inequality
with
$s = \nabla_{i_k} f(x_{k})$
and
$u = \frac{1}{a w_a(\delta (k+1))^{1/p}}(x_k-x^\ast)$.
Both descent conditions~\eqref{eq:prog1} imply $\|x_k-x^\ast\| \le R$, yielding the final inequality.
Therefore, we have shown that for all $k \ge 0$, $\mathbb{E}[E_{k+1}|x_k]-E_k \le c_p \delta R^p/a^p.$
This implies
$\mathbb{E}[E_k] \le E_0 +  c_p \delta k R^p/a^p.$
Therefore
\begin{talign*}
\mathbb{E}[f(x_k)] - f(x^\ast) \leq \frac{f(x_0) - f(x^\ast)}{(1+\delta k/(ap))^p} + c_p \frac{R^p}{a^p}\frac{\delta k}{(1+\delta k/(ap))^p}.
\end{talign*}
Since $a > 0$ was arbitrary, we may choose $a = R\frac{(c_p\delta k)^{1/p}}{(f(x_0) - f(x^\ast))^{1/p}}$ to obtain the bound
\begin{talign*} 
\mathbb{E}[f(x_k)] - f(x^\ast) 
	\leq \frac{2(f(x_0) - f(x^\ast))}{\left(1+\frac{(f(x_0) - f(x^\ast))^{1/p}}{Rc_p^{1/p}p}(\delta k)^{\frac{p-1}{p}}\right)^p}
	=\textstyle O(1/(1+\frac{1}{Rp}{(\delta k)^{\frac{p-1}{p}}})^{p})
\end{talign*}
as desired. 
\subsubsection{Proof of Theorem~\ref{thm:333}}
Take the energy function
$E_k = f(x_k) - f(x^\ast),$
and observe that if~\eqref{eq:p1} holds, then we have:
\begin{talign*} 
\frac{\mathbb{E}[E_{k+1}|x_k] - E_k}{\delta} = \frac{\mathbb{E}[f(x_{k+1})|x_k] - f(x_k)}{\delta}&\overset{\eqref{eq:corprog1}}{\leq} -\mathbb{E}[\|\nabla_{i_k} f(x_k)\|_{\ast}^{\frac{p}{p-1}}|x_k]\\ &= -\frac{1}{d}\sum_{i=1}^d\| \nabla_i f(x_k)\|^{\frac{p}{p-1}}\\
&\leq -\frac{1}{d}\| \nabla f(x_k)\|^{\frac{p}{p-1}}\\
& \overset{\eqref{eq:gdom}}{\leq} -\frac{1}{d}\frac{p}{p-1}\mu^{\frac{1}{p-1}} E_k,
\end{talign*} 
or rewritten, $\mathbb{E}[E_{k+1}] \leq \left(1 -\frac{1}{d}\frac{p}{p-1} \mu^{\frac{1}{p-1}}\delta\right)E_k$. Summing gives the bound \[\mathbb{E}[E_{k+1}] \leq \left(1-\frac{1}{d}\frac{p}{p-1}\mu^{\frac{1}{p-1}}\delta\right)^kE_0 \leq e^{-\frac{1}{d}\frac{p}{p-1}\mu^{\frac{1}{p-1}}\delta k}E_0,\] %
\subsubsection{Rescaled coordinate descent}
Rescaled coordinate descent, 
\begin{align}\label{eq:rgd-cor}
x_{k+1} = x_k - \eta_{i_k}^{\frac{1}{p-1}}\frac{\nabla_{i_k}f(x_k)}{\|\nabla_{i_k} f(x_k)\|^{\frac{p-2}{p-1}}} = \arg \min_{x \in \X} \left\{ \langle \nabla_{i_k} f(x_k), x\rangle + \frac{1}{\eta_{i_k} p}\|x - x_k \|^p\right\}
\end{align}
where $0<\eta_{i_k}<\infty $ for $ i_k \in \{1,\dots k\}$, satisfies~\eqref{eq:corprog1} provided the objective is strongly smooth along each coordinate direction.
\begin{definition}
A function $f$ is {\bf strongly smooth} of order $p$ along each coordinate direction for $p>1$, if there exist constants $0<L_1^{(i)}, \dots, L_p^{(i)}<\infty$ for $i = 1, \dots, d$, such that for $m = 1, \dots, p-1$ and for all $x \in \mathbb{R}^d$, as well as for all $i \in \{1, \dots d\}$ \begin{align}\label{eq:ss-cor}
\nabla^m f(x) (\nabla_{i} f(x))^m \leq L_m^{(i)} \|\nabla_{i} f(x)\|_\ast^{m + \frac{p-m}{p-1}}
,\end{align}
and moreover for $m = p$, $f$ satisfies the condition $\|\nabla^p f(x)\|\leq L_p^{(i)}$.
\end{definition} 
We summarize our results regarding the rescaled coordinate descent in the following Lemma.
\begin{lemma}\label{lem:rgd-cor}
Suppose $f$  is strongly smooth of order $p \geq 2$ along each coordinate direction with constants $0<L_1^{(i)}, \dots, L_p^{(i)} < \infty$ for $i = 1,\dots, d$. Then rescaled gradient descent~\eqref{eq:rgd-cor} with step size
\begin{align}\label{eq:step-2}
0 < \eta_i^{\frac{1}{p-1}} \leq \min\left\{1, \frac{1}{\left(2\sum_{m=2}^p \frac{L_m^{(i)}}{m!}\right)}\right\}
\end{align}
 satisfies~\eqref{eq:corprog1} with $\delta = \min_{i = 1,\dots,d}\eta_i^{\frac{1}{p-1}}/2$. %
\end{lemma}
\subsection{Accelerating Coordinate Descent Methods}
Coordinate descent algorithms of order $p$ can also be accelerated.%
Suppose $f$ is convex. 
Set $A_k = C\delta^p k^{(p)}$  where we use the rising factorial $ k^{(p)} = k(k+1)\cdots(k+p-1) $. Denote $\alpha_k := \frac{A_{k+1} - A_k}{\delta} = Cp\delta^{p-1}(k+1)^{(p-1)}$ and $\tau_k := \frac{\alpha_k}{A_{k+1}} = \frac{k}{\delta(k+p)}$. We write the algorithm as,
\begin{subequations}\label{eq:nestcor1}
\begin{talign}
x_{k} &= \delta\tau_k z_k + (1 - \delta\tau_k) y_k\label{eq:up1111}\\
z_{k+1} &= \arg \min_{z}\left\{ \alpha_k\langle \nabla_{i_k} f(x_k), z\rangle +  \frac{1}{\delta} D_h(z, z_k)\right\}\label{eq:up2222}
\end{talign}
\end{subequations} 
where the update for $y_{k+1}$ satisfies the descent condition 
\begin{talign}\label{eq:c11}
\frac{f(y_{k+1}) - f(x_k)}{\delta^{\frac{p}{p-1}}} \leq -  \|\nabla_{i_k} f(x_k)\|^{\frac{p}{p-1}}.
\end{talign}
For algorithm~\eqref{eq:nestcor1}, using~\eqref{eq:lyap} we compute
\begin{talign}\label{eq:pro1111}
 \frac{E_{k+1} - E_k}{\delta} = \frac{D_h(x^\ast, z_{k+1}) - D_h(x^\ast, z_k)}{\delta}+ \frac{A_{k+1}}{\delta}(f(y_{k+1}) - f(x^\ast)) - \frac{A_{k}}{\delta} (f(y_{k}) - f(x^\ast)).
\end{talign}
We bound the first part,
 \begin{talign}
\frac{D_h(x^\ast, z_{k+1}) - D_h(x^\ast, z_k)}{\delta} &= - \left\langle \frac{\nabla h(z_{k+1}) - \nabla h(z_k)}{\delta}, x^\ast - z_{k+1}\right\rangle - \frac{1}{\delta} D_h(z_{k+1}, z_k)\notag\\
&\overset{\eqref{eq:up2222}}{=}\alpha_k\langle \nabla_{i_k} f(x_{k}), x^\ast - z_k\rangle +\alpha_k \langle \nabla_{i_k} f(x_{k}), z_k - z_{k+1}\rangle\notag\\
&\quad - \frac{1}{\delta} D_h(z_{k+1}, z_k)\notag\\
& \leq \alpha_k\langle \nabla f(x_{k}), x^\ast - z_k\rangle - \xi_k -(\delta/m)^{\frac{1}{p-1}}\alpha_k^{\frac{p}{p-1}}\|\nabla_i f(x_{k})\|^{\frac{p}{p-1}}, \label{eq:pro2222} 
\end{talign}
where $\xi_k = \alpha_k\langle \nabla f(x_{k}) - \nabla_{i_k}f(x_k), x^\ast - z_k\rangle$ which is a martingale. The inequality follows from the $m$-uniform convexity of $h$ of order $p$ and the Fenchel-Young inequality $\langle s,u\rangle + \frac{1}{p}\|u\|^p \geq -\frac{p}{p-1} \|s\|_\ast^{\frac{p}{p-1}}$, with $u =(m/\delta)^{\frac{1}{p}}(z_{k+1} - z_k)$ and $s = (\delta/m)^{\frac{1}{p}}\alpha_k^{\frac{p}{p-1}} \nabla_{i_k} f(x_k)$. Plugging in update~\eqref{eq:up1}, 
\begin{talign}
 \alpha_k\langle \nabla f(x_{k}), x^\ast - z_k\rangle &= \alpha_k\langle \nabla f(x_{k}), x^\ast - y_k\rangle + \frac{A_{k+1}}{\delta}\langle \nabla f(x_{k}), y_k - x_{k}\rangle \notag\\
 &= \alpha_k\langle \nabla f(x_{k}), x^\ast - x_k\rangle + \frac{A_{k}}{\delta} \langle \nabla f(x_k), y_k - x_k\rangle\notag\\
& \leq - \left(\frac{A_{k+1}}{\delta}(f(y_{k+1}) - f(x^\ast)) -\frac{A_{k}}{\delta} (f(y_{k}) - f(x^\ast))\right) \notag\\
&\quad+  A_{k+1}\frac{f(y_{k+1}) - f(x_k)}{\delta}\notag\\
&\overset{\eqref{eq:c11}}{\leq} -  \left(\frac{A_{k+1}}{\delta}(f(y_{k+1}) - f(x^\ast)) -\frac{A_{k}}{\delta} (f(y_{k}) - f(x^\ast))\right) \notag\\
&\quad-  A_{k+1}\delta^{\frac{1}{p-1}}\|\nabla_{i_k} f(x_k)\|^{\frac{p}{p-1}}.\label{eq:pro3333}
\end{talign}
The first inequality follows from the convexity of $f$ and rearranging terms. The second inequality uses~\eqref{eq:c11}. Combining~\eqref{eq:pro1111} with \eqref{eq:pro2222} and \eqref{eq:pro3333} we have, 
\begin{talign*} 
\frac{E_{k+1} - E_k}{\delta} \leq \left((\delta/m)^{\frac{1}{p-1}}(Cp\delta^{p-1}(k+1)^{(p-1)})^{\frac{p}{p-1}}  -C\delta^{\frac{1}{p-1}}\delta^p(k+1)^{(p)}\right)\|\nabla_{i_k} f(x_{k})\|^{\frac{p}{p-1}}  - \xi_k.
\end{talign*}
Given $((k+1)^{(p-1)})^{\frac{p}{p-1}}/(k+1)^{(p)}\leq 1$, it suffices that
$C \leq 1/mp^{p} $ to ensure $\frac{\mathbb{E}[E_{k+1}|x_k] - E_k}{\delta} \leq 0$. Summing, we obtain the desired bound. 
\begin{align*} 
\mathbb{E}[f(x_k)] - f(x^\ast) \lesssim 1/(\delta k)^p.
\end{align*} 
\subsubsection{Accelerating rescaled coordinate descent}
\label{app:argd_cor}
A corollary to the coordinate descent property of rescaled descent with step size~\eqref{eq:step-2} is that it can be combined with sequences~\eqref{eq:up1111} and~\eqref{eq:up2222} to form a method with an $O(1/(\delta k)^p)$ convergence rate upper bound.
\begin{algorithm}[ht!]
\caption{Nesterov-style accelerated rescaled coordinate descent.}
\label{eq:argd_cor}
\begin{algorithmic}[1]
\Require{$f$ is {\em strongly smooth of order $p$} along each coordinate direction and $h$ satisfies $D_h(x,y) \geq \frac{1}{p}\|x - y\|^p$.}
	\State Set $x_0 = z_0 =0$ and $A_k = C\delta^p k^{(p)}$, $\alpha_k =  \frac{A_{k+1} - A_k}{\delta} =  Cp\delta^{p-1}(k+1)^{(p-1)} $ and $\tau_k = \frac{\alpha_k}{A_{k+1}} = \frac{k}{\delta(k+p)}$ where $ k^{(p)} := k(k+1)\cdots(k+p-1)$.\\
	 {\bf for} $k = 1, \dots, K$ {\bf do}
\State $x_{k} =  \delta \tau_k z_k + (1 - \delta \tau_k) y_k$
\State sample $i_k \in \{1, \dots, d\}$. Update
\State $z_{k+1} = \arg \min_{z}\left\{ \alpha_k\langle \nabla_{i_k} f(x_k), z\rangle +  \frac{1}{\delta}D_h(z, z_k)\right\}$ 
\State $y_{k+1} =x_k -  \eta_{i_k}^{\frac{1}{p-1}}\frac{\nabla_{i_k} f(x_k)}{\|\nabla_{i_k} f(x_k)\|_\ast^{\frac{p-2}{p-1}}} $
\State \Return $y_K$.
\end{algorithmic}
\end{algorithm}
We summarize this result in the following theorem.
\begin{theorem}
Suppose $f$ is convex and strongly smooth of order $1<p<\infty$ along each coordinate direction $i$ with constants $0<L_1^{(i)},\dots,L_p^{(i)}<\infty$. Also suppose $\eta_i$ satisfies~\eqref{eq:step-2}. Then Algorithm~\ref{eq:argd_cor} satisfies,
\begin{align*} 
\mathbb{E}[f(y_k)] - f(x^\ast) \lesssim   1/(\delta k)^{p}. 
\end{align*}  
\end{theorem}

\subsection{Optimal Universal Higher-order Tensor Methods}
\label{app:Uni}
We say that it has H\"older continuous $(p-1)$-st order gradients of degree
$\nu \in [0, 1]$ on a convex set $\mathcal{X} \subseteq \text{dom} f$, if for some constant $L_\nu$ it holds
\begin{talign}\label{eq:hold}
\|\nabla^{p-1} f(x) - \nabla^{p-1} f(y)\| \leq L_\nu\|x - y\|^\nu
\end{talign}

The final result of our paper contains the analysis of the following optimal algorithm for minimizing functions that satsify~\eqref{eq:hold}
\begin{algorithm}[H]
\caption{Monteiro-Svaiter-style universal higher-order tensor method.}
\label{alg:argd3}
\begin{algorithmic}[1]
\Require{$f$ satisfies~\eqref{eq:hold} with parameters $p$ and $L_\nu$, $h$ is $1$-strongly convex, $B = I, \tilde p = p-1 + \nu$.}
	\State Set $x_0 = z_0 =0$, $A_0 = 0$, $\delta^{\frac{3p-2}{2}} = \eta$,\,$\eta= L_\nu/(p-2)!$\\
	{\bf for} $k = 1, \dots, K$ {\bf do}
	\begin{subequations}
		\State Choose $\lambda_{k+1}$ (e.g. by line search) such that 
	\begin{talign}\label{eq:line-search}
	\frac{1}{2} \leq \frac{\lambda_{k+1} \|y_{k+1} -  x_k\|^{\tilde p-2}}{\eta} \leq \frac{3}{4},
	\end{talign}
	where 
	\begin{talign}\label{eq:hold_hgd} 
	y_{k+1} = \arg \min_{x\in \mathcal{X}} \left\{ f_{p-1}(x; x_k) + \frac{1}{\tilde{p}\eta }\|x - x_k\|^{\tilde{p}}\right\},
	\end{talign}
	\end{subequations}
	and $\alpha_{k} =\frac{ \lambda_{k+1} + \sqrt{\lambda_{k+1} + 4A_k\lambda_{k+1}}}{2\delta}$, $A_{k+1} = \delta \alpha_{k}  + A_k$, $\tau_k = \frac{\alpha_k}{A_{k+1}}$ (so that $\lambda_{k+1} = \frac{\delta^2 \alpha_k^2}{A_{k+1}}$) and 
		$$\textstyle x_k = \delta \tau_k  z_k + (1-\delta \tau_k) y_k.$$

	\State Update $z_{k+1} = \arg \min_{z\in \mathcal{X}}\left\{ \alpha_k \langle \nabla f(y_{k+1}), z\rangle +   \frac{1}{\delta}D_h(z, z_k)\right\}$
\State \Return $y_K$.
\end{algorithmic}
\end{algorithm}
We summarize results on performance of Algorithm~\ref{alg:argd3} in the following corollary to Theorem~\ref{prop:MS}:
\begin{theorem}\label{thm:M-AHGD} Assume $f$ is convex and has H\"older continuous $(p-1)$-st order gradients. Then Algorithm~\ref{alg:argd3} satisfies the convergence rate upper bound
\begin{align*} 
f(y_k) - f(x^\ast) =O\left(1/(\delta k)^{\frac{3(p-1+\nu)-2}{2}}\right).
\end{align*}
\end{theorem}

To prove Theorem~\ref{thm:M-AHGD}, the first thing to notice is that the proof of~\cref{prop:MS} holds for all $\mathbb{R} \ni p>0$. Subsequently, to extend our analysis to Algorithm~\eqref{alg:argd3}, it is sufficient to show (1)~\eqref{eq:hold_hgd}  with the line search step~\eqref{eq:line-search} satisfies 
\begin{align}\label{eq:first_step} \|y_{k+1} - x_k - \lambda_{k+1}\nabla f(y_{k+1})\| \leq \frac{1}{2}\|y_{k+1} - x_k \|
\end{align} and that (2) there exists a sequence $(\lambda_{k+1}, y_{k+1})$ that satisfies~\eqref{eq:hold_hgd} and~\eqref{eq:line-search} simultaneously. 

\paragraph{(1)} Observe that the optimality condition for~\eqref{eq:hold_hgd} satisfies 
\begin{align*}
\nabla f_{p-1}(y_{k+1}; x_k) - \frac{1}{\eta}(y_{k+1} - x_k) \|y_{k+1} -  x_k\|^{\tilde p-2} = 0.
\end{align*}
so that $\|\nabla f_{p-1}(y_{k+1}; x_k) \| = \frac{1}{\eta}\|y_{k+1} -  x_k\|^{\tilde p-1}$. In particular,
\begin{align*}
y_{k+1} - x_k + \lambda_{k+1}\nabla f(y_{k+1}) =  \lambda_{k+1}\nabla f(y_{k+1}) - \frac{\eta}{\|y_{k+1} -  x_k\|^{\tilde p-2}}\nabla f_{p-1}(y_{k+1}; x_k).
\end{align*}
From the integral form of the mean value theorem it follows that 
\begin{talign*}
\|\nabla f_{p-1}(y; x) - \nabla f(y)\| \leq \frac{L_\nu}{(p-2)!}\|y - x\|^{p-2+\nu}.
\end{talign*}
Subsequently
\begin{talign*}
\|y_{k+1} - x_k + \lambda_{k+1}\nabla f(y_{k+1}) \| &\leq \lambda_{k+1} \frac{L_\nu}{(p-2)!}\|y_{k+1} - x_k\|^{\tilde p - 1} + \left|\lambda_{k+1} -\frac{\eta}{\|y_{k+1} -  x_k\|^{\tilde p-2}}\right| \|\nabla f_{p-1}(y_{k+1}; x_k)\|\\
& \leq \|y_{k+1} -  x_k\|\left( \lambda_{k+1} \frac{L_\nu}{(p-2)!}\|y_{k+1} - x_k\|^{\tilde p - 2} + |\frac{\lambda_{k+1}}{\eta} \|y_{k+1} - x_k\|^{\tilde p - 2} + 1| \right)
\end{talign*}
If we choose $\eta = L_\nu/(p-2)!$ and plug in our line search criterion~\eqref{eq:line-search}, we see condition~\eqref{eq:first_step} is met. %

\paragraph{(2)} We now show there exists a pair $(\lambda_{k+1}, y_{k+1})$ that satisfies~\eqref{eq:hold_hgd} and~\eqref{eq:line-search} simultaneously. This claim follows directly form the argument given by Bubeck et al~\cite[Sec 3.2]{Bubeck}, which did not rely on $p>0$ being an integer. For self-containment, we reproduce the argument here.
\begin{lemma} Let $A\geq 0$, $x,y \in \mathbb{R}^d$ such that $f(x) \neq f(x^\ast)$. Define the following functions:
\begin{talign*}
a(\lambda) &= \frac{\lambda + \sqrt{\lambda^2 + 4\lambda A}}{2} \\
x(\lambda) & = \frac{a(\lambda)}{A + a(\lambda)} x + \frac{A}{A + a(\lambda)} y\\
y(z) &= \arg \min_{x\in \mathcal{X}} \left\{ f_{p-1}(w; z) + \frac{1}{\tilde{p}\eta }\|w - z\|^{\tilde{p}}\right\}\\
g(\lambda) &= \lambda\|y(x(\lambda)) - x(\lambda)\|^{\tilde p-1}.
\end{talign*} 
Then we have $g(\mathbb{R}_+) = \mathbb{R}_+$.
\end{lemma}
The first claim is that $g(\lambda)$ is a continuous function of $\lambda$ . This follows from the fact that $y(z)$ is a continuous function of $z$. Furthermore, $g(0) = 0$, and since $f(x) \neq f(x^\ast)$ we also have $y(x) \neq x$ which proves $g(+\infty)= +\infty$

\begin{remark} The same binary line search step introduced by~\citet[Sec 4]{Bubeck} finds a $\lambda_{k+1}$ satisfying~\eqref{eq:line-search}. The argument given there did not rely on the fact that $p \in \mathbb{Z}_+$. \end{remark} %

\end{document}